\documentclass[11pt,english]{article}
\usepackage[T1]{fontenc}
\usepackage[latin9]{inputenc}
\usepackage[a4paper]{geometry}
\usepackage{xcolor}
\usepackage{pdfcolmk}
\usepackage{array}
\usepackage{float}
\usepackage{multirow}
\usepackage{amsmath}
\usepackage{amsthm}
\usepackage{graphicx}
\usepackage{setspace}
\usepackage[numbers]{natbib}
\PassOptionsToPackage{normalem}{ulem}
\usepackage{ulem}
\makeatletter

\providecommand{\tabularnewline}{\\}
\providecolor{lyxadded}{rgb}{0,0,1}
\providecolor{lyxdeleted}{rgb}{1,0,0}

\DeclareRobustCommand{\lyxsout}[1]{\ifx\\#1\else\sout{#1}\fi}

\newcommand{\lyxaddress}[1]{
	\par {\raggedright #1
	\vspace{1.4em}
	\noindent\par}
}

\date{}
\usepackage{algpseudocode}
\usepackage{amsmath}
\usepackage{lscape}
\usepackage{newunicodechar}
\newunicodechar{'}

\@ifundefined{showcaptionsetup}{}{%
 \PassOptionsToPackage{caption=false}{subfig}}
\usepackage{subfig}
\makeatother

\usepackage{babel}
\begin{document}
\title{\textbf{Drone Deployment Optimization for Direct Delivery with Time
Windows and Battery Replacements}}
\author{Tanveer Hossain Bhuiyan{*} $^{a,b}$, Mohammad Roni, $^{a}$ and Victor
Walker$^{a}$}
\maketitle

\lyxaddress{\begin{center}
$^{a}$ Energy and Environment Science \& Technology Department,
Idaho National Laboratory, Idaho, USA\vspace*{-5mm}
\par\end{center}}

\lyxaddress{\begin{center}
$^{b}$ Department of Mechanical Engineering, The University of Texas
at San Antonio, TX, USA\vspace*{-5mm}
\par\end{center}}

\lyxaddress{\begin{center}
{*}{\small{}Corresponding author email: tanveer.bhuiyan@utsa.edu}
\par\end{center}}
\begin{abstract}
Aerial drones offer a distinct potential to reduce the delivery time
and energy consumption for the delivery of time-sensitive and small
products. However, there is still a need in the relevant industry
to understand the performance of drone-based delivery under different
business needs and drone operating conditions. We studied a drone
deployment optimization problem for direct delivery of goods to customers
maintaining a specified time window. This paper presents a new mathematical
optimization-based decision-making methodology to help business owners
optimally route their drone fleet minimizing the total energy consumption,
required fleet size, and the required number of additional batteries.
A realistic feature of the optimization method is that instead of
replacing the drone battery after each return to the depot, it keeps
track of the remaining energy in the drone battery and decides on
battery replacements accounting for the drone routing and the user-specified
minimum required battery energy. Numerical results based on real drone
flight tests and delivery data provide insights into the effect of
different drone operating parameters on the energy consumption, required
fleet size, and the required number of battery replacements. Results
from a case study show that the total energy consumption, required
fleet size, and the required number of battery replacements increase
by 72.22\%, 22.2\%, and 200\%, respectively, as the drones fly over
the road networks compared to flying in a straight path. Additionally,
results show that using a mixed fleet of hexacopter and quadcopter
drones reduces the total energy consumption by 48.52\% compared to
using a homogeneous fleet of only hexacopters.
\end{abstract}

\section{Introduction}

This paper studies a drone deployment optimization problem for direct
delivery of goods from a retail business location (e.g., store) to
customer locations. The objective of the decision-maker (i.e., business
owner) is to find the optimal routing of drones that minimizes the
investment (e.g., drone fleet and battery costs) and operating (e.g.,
energy, labor, and maintenance costs) costs in delivering customer
orders within a specified time window. Specifically, the goal of this
study is to (1) develop a mathematical optimization-based decision-making
model that retail business owners can use to optimally deploy a fleet
of drones in delivering customer orders with a minimum required fleet
size (i.e., number of drones), required number of additional batteries,
and total energy consumption; (2) understand the effect of different
business and drone operating conditions on total energy consumption,
required fleet size, and required number of additional batteries;
and (3) provide numerical results and managerial insights based on
real drone flight test data to help business owners better understand
the potential of drones in goods delivery and evaluate alternative
decisions.

\subsection{Motivation}

Unmanned aerial vehicles (UAVs), also known as drones, have gained
abundant attention from businesses and consumers. For instance, Amazon
and several other companies (e.g., Google, DHL) have developed their
own drone delivery systems \citep{cheng2020drone}. Drones possess
several advantages compared to other vehicles, such as low capital
and operating costs, the capability of taking off and landing in a
limited space without requiring expensive infrastructure, and eliminating
the necessity of an operator on board \citep{gentili2022locating}.
Due to these advantages, drones offer a great benefit in numerous
transportation and logistics activities including but not limited
to traffic surveillance \citep{chow2016dynamic}, consumer product
delivery \citep{murray2015flying}, emergency medical supplies \citep{gentili2022locating},
and agriculture and forest fire surveillance \citep{clarke2014understanding}.

Specifically, drones offer a distinct potential to improve the delivery
of goods in the last-mile, especially for time-sensitive, small, and
localized deliveries. In recent years, general e-commerce and retail
online deliveries are growing, and in some areas they are approaching
a critical stage, such as in New York City, where over 1.5 million
daily deliveries are putting a strain on the road networks and curb
resources \citep{NY_times}. In addition, local business-to-consumer
deliveries are growing quickly for time-sensitive products, such as
prepared food, medicine, and groceries. With the growing online retail
business, cost-efficient and timely delivery of products from stores
to the customers\textquoteright{} homes is becoming a challenge for
business owners, especially in a densely populated urban environment
due to congestion in the road networks. Moreover, current delivery
services are also expensive, costing restaurants up to 30\% of each
order\textquoteright s value for food delivery \citep{wingReport}.
In these circumstances, drone-based delivery can significantly reduce
traffic congestion by avoiding the use of road networks, as well as
increase the delivery range, market reach, and sales volume for the
retail business owners by reducing the cost and delivery time of the
current system of using ground vehicles. Moreover, McKinsey has suggested
that 80\% of last-mile deliveries will be done by autonomous vehicles
within 8 years \citep{McKinsey}, and that drones will likely be the
chosen method for a very large segment of the population, especially
for suburban and rural consumers needing same-day and speedy deliveries.

Despite the great potential and benefit, drone-based delivery suffers
from a major limitation\textemdash limited battery capacity\textemdash that
limits the drone delivery range. Therefore, in a drone routing problem,
battery capacity, energy consumption, and battery replacement decisions
(i.e., swapping the depleted battery with a fully-charged one) need
to be incorporated to ensure safe deliveries of products. Most research
studies on drone deployment and operations assume a fixed flight duration
to reflect the battery capacity \citep{cheng2020drone}, whereas some
studies consider the drone energy consumption to be a function of
only drone speed, or package weight. But, in practice, drone energy
consumption, and thus, the battery capacity is affected by several
drone operating factors\textemdash speed, package weight, battery
weight, distance traveled, and duration of flight components (e.g.,
ascend, descend, hover). Therefore, all these factors need to be accounted
for while computing the drone energy consumption and in the drone
routing decisions. Also, depending on the type of products, delivery
locations, and customer preferences, business owners may need to use
different delivery methods, such as landing with the package on the
ground, or hovering above the ground and delivering the product through
a winch. Different delivery methods cause different amounts of energy
consumption, which can also affect the drone range. Therefore, practitioners
(e.g., business owners) need to get insight into the question: how
much do the package weight, drone speed, and delivery method affect
the drone energy consumption, required fleet size (i.e., number of
drones), and the required number of battery replacements?

In the drone routing literature on direct delivery from business to
customers, all the studies assume that drones always fly in a straight
path. However, in some areas, especially in densely populated urban
areas, drones may not be able to fly in a straight path due to the
presence of no-fly zones in the flight path, such as schools, large
buildings, and restricted private property. To avoid these no-fly
zones, flying over the road networks is a viable option as they are
public property, resulting in a higher amount of energy consumption
and delivery time. Therefore, how does the drone flight path (e.g.,
straight vs. over the road networks) affect the drone delivery range,
energy consumption, required fleet size, and the required number of
battery replacements in delivering all customer orders?

Also, all existing studies on direct delivery from business to customer
consider that each drone delivers packages to multiple customers in
a single trip before returning to the depot, which is not always possible,
especially for rotary drones, such as quadcopters and hexacopters.
Moreover, in current drone-delivery systems, such as Amazon Prime
Air press releases \citep{AmazonPrime} demonstrate drones delivering
packages to a single customer on each trip. In case of delivering
a single customer order on each trip, replacing the drone battery
after each return to the depot is not practical. Because it is likely
that there is a sufficient amount of energy remaining in the drone
battery to make an additional delivery to another customer. Therefore,
replacing the battery after each return can lead to an unnecessarily
large number of battery replacements and less drone utilization, causing
the business owner to have a large number of additional batteries
and battery charging infrastructure, as well as a larger number of
drones in the fleet to maintain a desired delivery time window for
customers. Therefore, the question remains on how the drones be routed
to minimize the required fleet size? How should drones be routed to
minimize the required number of battery replacements?

Most studies assume a homogeneous fleet of drones, where all drones
are identical. But, customers are usually located at different distances
from the business location (e.g., store) and their orders have different
weights. Therefore, different drones may be suitable for different
customer orders in terms of energy-efficiency. Therefore, how much
is the benefit of using a mixed (heterogeneous) fleet of drones compared
to using a homogeneous fleet?

The above-mentioned questions need to be answered quantitatively using
real data for the retail business owners to demonstrate the potential
of drones in the direct delivery of products to customers. Motivated
by these needs, we study a drone routing problem for direct delivery
of goods and proposed a mathematical optimization-based approach that
incorporates realistic drone operating conditions and provides insights
to the questions based on real drone flight test data. Specifically,
our proposed mathematical optimization approach can provide insights
into the effect of different operating parameters\textemdash drone
speed, package weight, flight path, delivery method, fleet type\textemdash on
the delivery range, required fleet size, required number of battery
replacements, and energy consumption.

\subsection{Related Literature}

Drone routing optimization has been gaining a lot of attention to
researchers recently. Several research teams have conducted comprehensive
review of the recent studies on drone delivery and operations. Rojas
et al. \citep{rojas2021unmanned} and Marcina et al. \citep{macrina2020drone}
provide a review of the research on generic drone routing problems
and for parcel delivery, respectively. Boysen et al. \citep{boysen2021last}
provides a review of the operations research methods used for solving
the last-mile delivery problems with drones. Otto et al. \citep{otto2018optimization},
Coutinho et al. \citep{coutinho2018unmanned}, and Chung et al. \citep{chung2020optimization}
provide comprehensive reviews of the optimization approaches\textemdash models
and algorithms\textemdash used for drone routing problems, where Chung
et al. \citep{chung2020optimization} focuses on drone and drone-truck
combined routing applications, and Otto et al. \citep{otto2018optimization}
focuses on the civilian applications (e.g., medical supplies and disaster
management) of drones. Meanwhile, Barmpounakis et al. \citep{barmpounakis2016unmanned}
provides a review of the studies on using drones for surveillance
in transportation and traffic engineering.

In the drone-based product delivery literature, combined drone and
ground vehicle routing problems have been mostly studied, where ground
vehicles (e.g., trucks, delivery van) carry one or multiple drones
to perform localized deliveries. Researchers have studied a different
variants of this combined drone-truck delivery problem: (1) drone
and truck both deliver packages to customers \citep{poikonen2017vehicle,wang2017vehicle};
and (2) packages are delivered by drone only, whereas the trucks act
as a mobile depot to carry drones to a longer distance \citep{mathew2015planning}
and provide drone battery recharging/replacement services \citep{tokekar2016sensor}.
Wang et al. \citep{wang2017vehicle} were the first to model a combined
drone-truck parcel delivery problem seeking to reduce the total delivery
time, where drones are carried by trucks. This proposed model ignores
several realistic features of the problems, including the limitation
of drone battery life, cost, and different distance metrics for drones
and trucks. However, this study presents a comparison of the delivery
time between a truck-only fleet and a drone-truck combined fleet.
The Wang et al. \citep{wang2017vehicle} study was extended by Poikonen
et al. \citep{poikonen2017vehicle} by incorporating the limited drone
battery life and distance metrics for trucks and drones. Di Puglia
Pugliese and Guerriero \citep{di2017last} extended the work of Wang
et al. \citep{wang2017vehicle} and Poikonen et al. \citep{poikonen2017vehicle}
by incorporating delivery time windows for customers in a parcel delivery
problem aiming to minimize the overall travel cost. These authors
made numerical comparisons between the use of a combined drone-truck
fleet and a truck-only fleet in terms of transportation cost and total
delivery time for last-mile delivery.

Murray and Chu \citep{murray2015flying} were the first to introduce
two different variations of the drone-truck delivery problem, where
in one variant, truck and drones work synchronously in delivering
the package, whereas in the second variant, both the truck and drones
independently deliver packages to the customers. However, in this
routing problem, the authors only considered a single truck with multiple
drones. Ham \citep{ham2018integrated} and Murray and Raj \citep{murray2020multiple}
extended the drone-truck routing model of Murray and Chu \citep{murray2015flying}
by introducing multiple drones and trucks. These studies were further
extended in \citep{kitjacharoenchai2019multiple}, where the drones
can fly from a delivery truck and land on another truck, adding more
flexibility to the delivery problem. Heimfarth et al. \citep{heimfarth2022mixed}
extended the independent drone truck delivery problem accounting for
the customer preferences in delivering the parcel either using a truck
or a drone. The authors considered a less strict delivery time window
requirement for each delivery, where the customer is provided a rebate
if the time window is violated. There also exist studies where ground
vehicles are used only as a carrier for drones (e.g., \citep{huang2020round,gentili2022locating}).
Huang et al. \citep{huang2020round} proposed a drone routing model
for parcel delivery utilizing public transportation to carry drones
to a longer distance.

However, none of the above-mentioned studies modeled drone battery
replacements or charging phenomena, rather they made a limiting assumption
that a drone battery is enough to complete each trip. Also, these
studies mostly considered identical drones in the fleet. However,
using a mixed fleet of drones is beneficial in making energy-efficient
delivery. Gentili et al. \citep{gentili2022locating} studied the
problem of locating mobile depots (e.g., ground vehicles) and routing
drones for providing emergency medical supplies to disaster-affected
areas aiming to minimize the total disutility of the demand locations.
The authors considered a mixed fleet, consisting of two types of drones\textemdash short-range
and long-range. In their routing problem, each drone carries a single
package to the demand point and returns back to the mobile depot.
The authors did not explicitly model the battery replacement and resuming
of the drone energy level, rather just added a fixed amount of time
for package loading and battery replacement each time a drone returned
to the depot. However, replacing the drone battery at each return
is not realistic as there may be energy remaining in the drone battery
to make the next trip. Therefore, this solution can lead to a larger
number of unnecessary battery replacements, which results in a larger
number of additional batteries and battery charging infrastructure.
In the drone-truck hybrid delivery literature, none of the above studies
explicitly modeled the drone battery replacement to minimize the required
number of battery replacements, neither did they evaluate the impact
of no-fly zones on fleet sizing decisions, delivery time, and energy
consumptions. Only Jeong et al. \citep{jeong2019truck} incorporated
no-fly zones in the drone flight path in a drone-truck hybrid parcel
delivery problem and evaluated the effect of no-fly zones on delivery
time. However, these authors did not model the battery replacements
and mixed fleet of drones.

Compared to truck-drone hybrid delivery, relatively less work has
been done on deploying drones for the direct delivery of products
to customers. Yadav and Narasimhamurthy \citep{yadav2017heuristics}
studied a UAV routing problem seeking to minimize the total time for
delivery of all customer orders, where multiple drones can be used
to deliver the demand of a single customer. However, the authors made
an impractical assumption that drones can fly indefinitely in delivering
customer orders. But, in practice, battery energy capacity is the
major limitation of the drone-based deliveries that limits the drone
delivery range. A Genetic algorithm-based procedure was presented
by San et al. \citep{san2016delivery} for assigning drones to deliver
packages to customer locations. Later, Song et al. \citep{song2018persistent}
modeled the UAV delivery problem as a mixed-integer program and solved
using a heuristics algorithm. Unlike Yadav and Narasimhamurthy \citep{yadav2017heuristics}
and San et al. \citep{san2016delivery}, the authors considered the
effect of package weight on the flight time of UAVs to reflect the
limited battery capacity. However, they did not account for the other
factors (e.g., drone speed, flight path) in computing the energy consumption,
neither did they explicitly model battery energy consumption or battery
replacements.

Dorling et al. \citep{dorling2016vehicle} modeled the battery energy
consumption as a function of drone battery and package weight. The
authors proposed mixed-integer linear programming formulations for
two variants of a drone routing problem with the objective functions
of minimizing the total delivery time and the total cost, respectively.
The model assumes that each drone can deliver multiple packages to
different customers in a single trip and ensures that the battery
capacity is enough to make a trip before the trip starts. The drone
battery is replaced after each trip. However, this study does not
explicitly model the minimum energy requirement to determine battery
replacement, neither does it seek to minimize the number of battery
replacements. Rabta et al. \citep{rabta2018drone} modeled a drone
routing problem in addition to drone battery charging location decisions
for delivering disaster relief packages as a mixed-integer program.
Like Dorling et al. \citep{dorling2016vehicle}, each drone can deliver
multiple packages in a single trip and battery consumption is a function
of package weight. Their model tracks the remaining drone battery
energy in the model using a variable to determine the charging need.
But the battery charging was not modeled explicitly accounting for
the charging time or rate.

Some studies incorporated a time window for customer deliveries in
the drone routing problem and optimize the fleet size in delivering
customer orders. Troudi et al. \citep{troudi2018sizing} modeled a
capacitated vehicle routing problem with time windows as a mixed-integer
program to minimize the required number of drones and batteries. Similar
to Dorling et al. \citep{dorling2016vehicle}, the authors assumed
that each drone carries multiple packages in a single trip, and that
the battery is replaced after each trip. Exploiting these assumptions,
the authors minimized the number of batteries used by minimizing the
number of trips. Cheng et al. \citep{cheng2020drone} studied another
multi-trip (i.e., each drone performs multiple trips starting and
ending at the depot) drone routing problem with time windows aiming
to minimize the travel cost and electricity consumption cost related
to the drone battery. The authors considered drone energy consumption
as a function of package weight and distance traveled. However, they
used theoretical power consumption during hovering to approximate
the forward flight energy consumption, which is not realistic. Another
multi-trip drone delivery problem with time windows was studied by
Kong et al. \citep{kong2022trajectory} to minimize the total routing
distance in delivering all parcels to customers. Considering energy
consumption as a function of package weight, Choi and Schonfeld \citep{choi2017optimization}
studied an automated drone delivery system to minimize the total cost
by optimizing fleet size. The authors demonstrated the sensitivity
of the drone speed and battery capacity on the system cost and fleet
size, respectively. All the above-mentioned studies assume the drone
fleet is homogeneous.

Coelho et al. \citep{coelho2017multi} studied a drone routing problem
for a mixed fleet of drones, containing six different types of drones,
with the objectives of minimizing the total distance, total delivery
time, and the number of drones used. However, the authors considered
drone energy consumption as a function of drone speed only. But the
drone energy consumption is affected by several other operating parameters\textemdash package
weight, drone body and battery weight, distance traveled, and duration
of different flight segments (e.g., ascend, descend, hover). In addition
to Coelho et al. \citep{coelho2017multi}, all the above-mentioned
studies on drone-only delivery modeled the drone routing problem assuming
each drone delivers packages to multiple customers in a single trip,
and that the battery is replaced after each trip. Under this assumption,
the routing problem is modeled just to ensure the battery energy is
enough to complete each trip in visiting the customers. But, in practice,
many drones\textemdash especially rotary drones\textemdash are designed
to carry a single package on a trip. In this case when a drone visits
only a single customer on each trip, replacing the drone battery after
each return to the depot will lead to an unnecessarily large number
of battery replacements and less drone utilization, which can eventually
result in a large number of additional batteries or charging infrastructure
and the larger number of drones to maintain the same delivery time
window. Additionally, all the above studies on drone routing assume
drones fly in a straight path, which may not be always possible, especially
in a dense urban environment.

There exist other studies on drone routing focused on the effect of
drone delivery on carbon emission (e.g., \citep{chiang2019impact,figliozzi2017lifecycle})
and locating drone deployment and charging facilities to improve demand
coverage (e.g., \citep{chauhan2019maximum,trotta2018joint}).

\subsection{Contributions}

In summary, existing studies have not simultaneously considered the
different factors (e.g., drone speed, package weight, distance traveled,
duration of ascend, descend, hover) affecting drone energy consumption
and delivery range. A lack of studies exist analyzing the improvement
in the energy-efficiency of the drone fleet due to using a mixed fleet
of drones compared to a homogeneous fleet. In the direct delivery
of goods from business to customers, no research has explicitly modeled
the drone battery replacement decisions accounting for the remaining
energy in drone battery, user-specified minimum required energy in
drone battery, and drone routing when each drone visits a single customer
on each trip. Existing studies assume that each drone visits multiple
customers in a single trip before returning to the depot and replaces
the drone battery after each return to the depot, which can lead to
an unnecessarily large number of battery replacements when each drone
visits a single customer in each trip, eventually necessitating a
large number of additional batteries and drones. Additionally, no
research has analyzed the effect of routing drones over the road networks\textemdash due
to avoiding restricted and no-fly zones\textemdash on total energy
consumption, required fleet size (i.e., number of drones), and the
required number of battery replacements.

Therefore, to fill the gaps in the literature, we studied a drone
deployment optimization problem for direct delivery of goods that
accounts for the different factors (e.g., drone speed, package weight,
distance traveled, duration of ascend, descend, hover) affecting drone
energy consumption, pickup time window for customer orders, and battery
replacement decisions to evaluate the performance of drone-based delivery
under different drone operating parameters\textemdash speed, flight
path (e.g., straight vs. over the road networks), fleet type (e.g.,
homogeneous vs. mixed fleet), and package delivery method (e.g., landing
vs. package dropping). This is the first study to (1) model the drone
routing problem where the drone battery is not replaced after each
return to the depot, rather a battery replacement is decided by the
optimization model accounting for the remaining battery energy, minimum
required energy for safe operation, and drone routing to minimize
the number of battery replacements needed; and (2) evaluate the effect
of several key drone operating parameters\textemdash flight path,
fleet type, and minimum required energy in the drone battery\textemdash on
the delivery range, energy consumption, required fleet size, and the
number of times battery needs to be replaced.

Our paper extends the drone routing literature by introducing a new
mathematical optimization approach that incorporates real-life drone
operating characteristics and presents new insights based on real
drone flight test data. Specifically, in this paper, we made the following
contributions: (1) developed new mathematical optimization models
that incorporate the business limitations and drone operating conditions
and seek to optimally route the drones to minimize the total energy
consumption, required fleet size, and the required number of battery
replacements in delivering a set of customer orders; (2) proposed
two valid inequalities to improve the computational efficiency of
the mathematical model; and (3) provided numerical results and new
insights based on real data into the effect of drone speed, package
weight, flight path, minimum required energy, delivery method, fleet
type, and pickup time window on the drone delivery range, total energy
consumption, required fleet size, and the required number of battery
replacements.

\section{\label{sec:Problem-description}Problem Description}

\begin{figure}
\begin{centering}
\includegraphics[width=12cm,height=10cm,keepaspectratio]{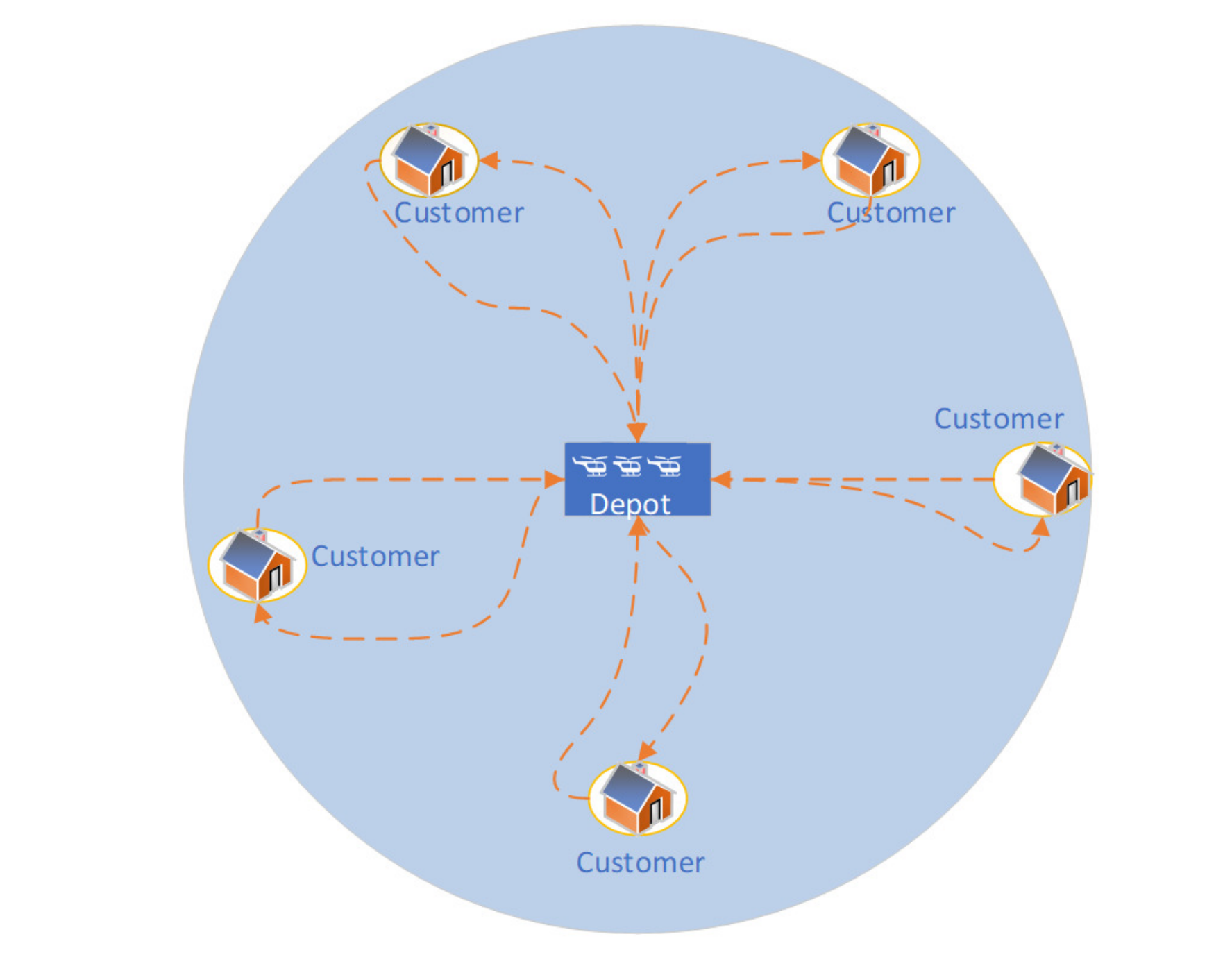}
\par\end{centering}
\caption{Business scenario.\label{fig:Business-scenario}}

\end{figure}

In this paper, we study the problem of a retail business owner (e.g.,
grocery stores, restaurants) seeking to efficiently deploy drones
to directly deliver goods to a set of customer locations that minimizes
the total energy consumption, the required fleet size, and the required
number of battery replacements. We assume that the business owner
owns the drone fleet and that all drones are stationed at the business
location that acts as a central depot for the drones. The drone fleet
may consists of the same type (i.e., homogeneous fleet) or different
types (i.e., mixed fleet) of drones. The customers are geographically
located around the central depot. The business owner knows the geographic
location (i.e., latitude and longitude) from where each customer order
originates, as well as the time stamp of order placement. In this
business scenario (see Figure \ref{fig:Business-scenario}), each
drone carries a single package in a flight from the depot to a delivery
(i.e., customer) location and then returns to the depot before flying
to the next delivery location.

Depending on the business and product types, customers may have different
preferences on the time window (e.g., 2 minutes for Starbucks coffee
vs. 30 minutes for a grocery item) they are willing to allow for drones
to pick up the products once they are ready for delivery. To account
for this realistic customer preference criterion in retail business,
we consider a pickup time window for each customer order within which
a drone must pick up the product from the depot. This time window
for each customer order is defined by the package ready time (i.e.,
earliest possible pickup time) and the maximum permissible delayed
pickup time. The travel time of a drone from the depot to a delivery
location consists of the time required to ascend to the flying height,
forward flight time to reach the delivery location, short hovering
duration required to ensure the correct landing location, and the
time required to descend to the ground. We also consider the time
required to load and unload the packages to/from the drone at the
depot and the delivery locations, respectively, as these package loading
and unloading times affect drone operation and satisfying the customer
order pickup time windows.

We consider that each drone starts with a fully-charged battery at
the beginning of the planning horizon. As the drones continue to deliver
packages to the customer locations, energy in the drone battery gets
consumed, and thus, the remaining energy in the drone battery gradually
decreases. The energy consumption during drone operation depends on
several operating parameters, such as drone speed, package weight,
distance traveled, duration of flight segments (i.e., ascend, descend,
hover), and delivery method. To ensure safe operation of the drones,
prior to leaving the depot for the next delivery location, we need
to ensure that the drone battery would have at least a minimum required
amount of energy remaining after returning from that delivery location.
Therefore, before leaving the depot for the next delivery location,
we need to compute the potential remaining energy in the drone battery
after returning from the next location. If the potential remaining
energy is less than the minimum required energy, the current battery
needs to be replaced with one that is fully-charged before leaving
for the next delivery location.

However, despite starting with a fully-charged battery to deliver
a package to a delivery location, if the potential remaining energy
in the drone battery drops below the minimum required energy, then
it is not possible for the drones to deliver a package to that particular
location. We refer this particular delivery location to be outside
the \textit{delivery range} of the drone. This \textit{delivery range}
of a drone varies with different operating parameters\textemdash drone
speed, package weight, drone flight path, delivery method, and minimum
required energy. A drone cannot be assigned to serve a location that
is outside its delivery range. In addition, each drone has a maximum
weight carrying capacity that cannot be exceeded, meaning that the
drone cannot be used to deliver a package that weighs more than the
maximum allowable weight.

\section{\label{sec:Model-Formulations}Methods}

To solve the drone deployment problem discussed in Section \ref{sec:Problem-description},
we formulate the problem as a mixed-integer program. In this section,
we present the mathematical formulations of the drone deployment problem
for both homogeneous and mixed fleets of drones. We lists the necessary
sets, parameters, and variables that support the mathematical formulations
in Table \ref{tab:Notation}.

\begin{table}[H]
\caption{Notation.\label{tab:Notation}}
\subfloat[\label{tab:Sets}Sets]{
\centering{}%
\begin{tabular}{cl}
\hline 
Sets & Description\tabularnewline
\hline 
\hline 
$\mathcal{I}$ & Set of delivery locations $i$\tabularnewline
$\mathcal{I^{\prime}}$ & Set of delivery locations including the dummy source $0$ and sink
$S$\tabularnewline
$\mathcal{D}$ & Set of drone types $d$\tabularnewline
\hline 
\end{tabular}}
\begin{centering}
\subfloat[Parameters\label{tab:Parameters}]{
\centering{}%
\begin{tabular}{c>{\raggedright}p{10cm}}
\hline 
Parameters & Description\tabularnewline
\hline 
\hline 
$C_{d}$ & Amortized cost of the drone type $d$\tabularnewline
$C_{bat}^{d}$ & Amortized cost of the battery of drone type $d$\tabularnewline
$v$ & Drone speed\tabularnewline
$C_{E}$ & Cost of per unit of energy\tabularnewline
$C_{L}$ & Wage of a drone operator\tabularnewline
$n_{d}$ & Number of drones of type $d$ a drone operator can operate simultaneously\tabularnewline
$C_{d}^{M}$ & Maintenance cost of a drone of type $d$\tabularnewline
$L_{i}$ & Distance of delivery location $i$ from depot\tabularnewline
$t_{Hi}=\frac{L_{i}}{v}$ & Time required for drone to arrive at delivery location $i$ from depot\tabularnewline
$t_{iH}=\frac{L_{i}}{v}$ & Time required for drone to return from delivery location $i$ to depot\tabularnewline
$t_{d}^{bat}$ & Time required to replace the battery of drone type $d$\tabularnewline
$e_{i}$ & Earliest possible pickup time\tabularnewline
$\Delta$ & Maximum permissible delay\tabularnewline
$\mathit{ch_{d}^{0}}$ & Initial energy in the battery of drone type $d$\tabularnewline
$\mathit{ch_{d}^{min}}$ & Minimum remaining energy required in the battery of drone type $d$\tabularnewline
$E_{dHi}$ & Energy consumption by a drone of type $d$ to arrive at delivery location
$i$ from depot with a package \tabularnewline
$E_{diH}$ & Energy consumption by drone of type $d$ to return empty from delivery
location $i$ to the depot\tabularnewline
$c_{ij}^{d}$ & Energy consumption by a drone of type $d$ to fly from location $i$
to $j$, i.e., $c_{ij}^{d}=E_{diH}+E_{dHj}$ \tabularnewline
$M_{d}^{max}$ & Maximum package weight carrying capacity of drone type $d$\tabularnewline
$M_{i}$ & Package weight for delivery location $i$\tabularnewline
\hline 
\end{tabular}}
\par\end{centering}
\centering{}\subfloat[Variables\label{tab:Variables}]{
\centering{}%
\begin{tabular}{c>{\raggedright}p{10cm}}
\hline 
Variables & Description\tabularnewline
\hline 
\hline 
$z_{ij}$ & 1 if delivery location $i$ is served immediately before $j$ by a
drone, 0 otherwise\tabularnewline
$y_{i}$ & 1 if drone battery is replaced after returning from delivery location
$i$, 0 otherwise\tabularnewline
$g_{i}$ & Remaining battery energy of a drone after returning from delivery
location $i$\tabularnewline
$g_{i}^{\prime}$ & Auxiliary variable storing the remaining battery energy after returning
from delivery location $i$\tabularnewline
$f_{i}$ & Timing of when a drone picks-up the package for delivery location
$i$ at depot \tabularnewline
$x_{di}$ & 1 if drone type $d$ is assigned to deliver a package to location
$i$, 0 otherwise \tabularnewline
\hline 
\end{tabular}}
\end{table}

\subsection{\label{subsec:Model_identical_drones}Model for Homogeneous Fleet
of Drones}

\begin{figure}[h]
\begin{centering}
\includegraphics[width=12cm,height=12cm,keepaspectratio]{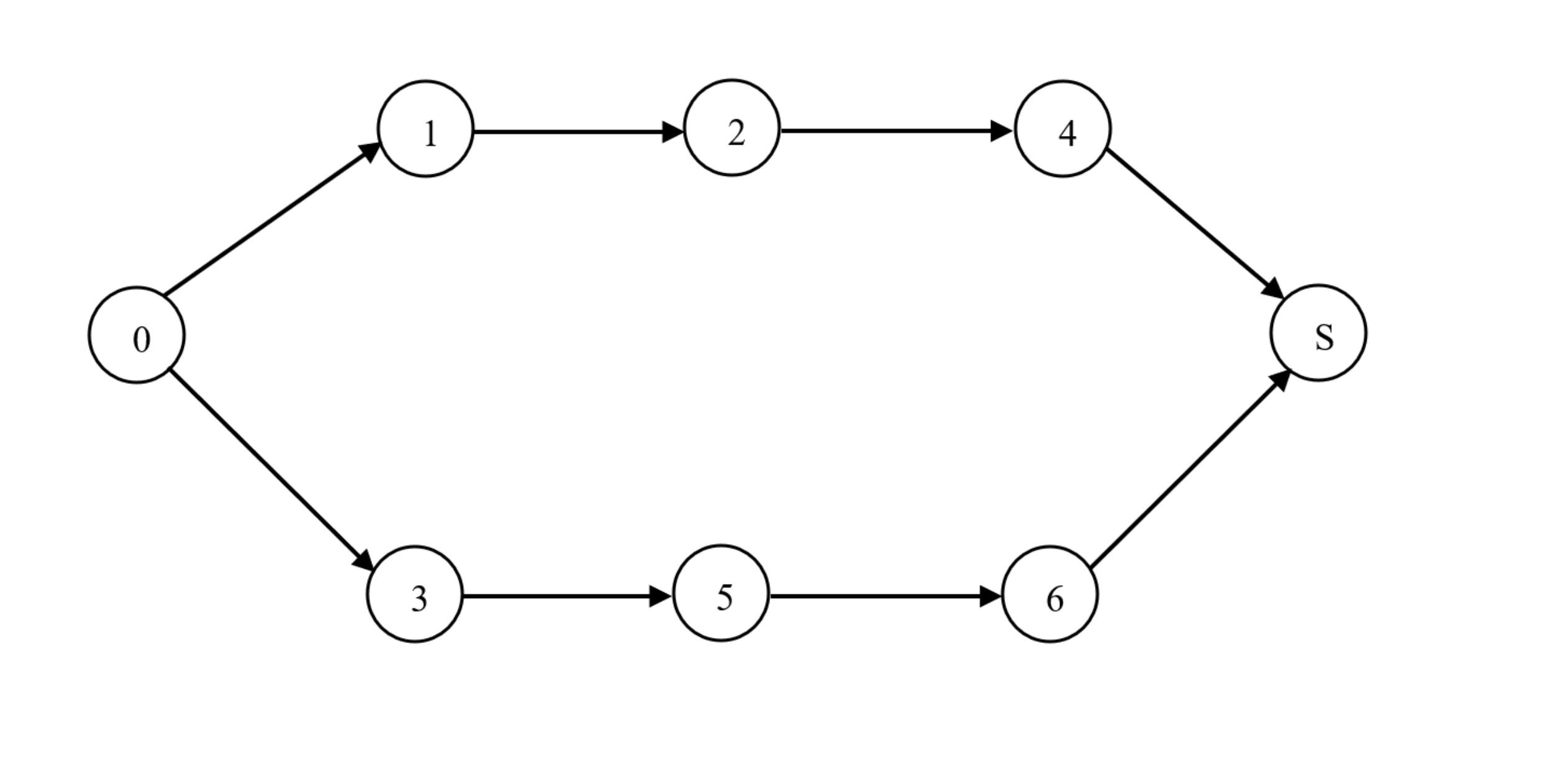}
\par\end{centering}
\caption{\label{fig:An-example-network_representation}An example network representation
of the original problem.}
\end{figure}

In this subsection, we assume that the drone fleet is homogeneous,
meaning that the drones are identical. Taking advantage of this assumption,
we can develop the drone deployment model as a two-index formulation,
similar to the vehicle routing literature (e.g., \citep{furtado2017pickup}).
As mentioned in Section \ref{sec:Problem-description}, in the retail
business scenario studied in this paper, each drone directly flies
from the depot to a delivery location and returns to the depot before
flying to the next location. To model this drone routing problem as
a network optimization problem (an example network is shown in Figure
\ref{fig:An-example-network_representation}), we add a dummy source
($0$) and a sink ($S$) to the set of delivery locations $\mathcal{I}$,
resulting in a new set $\mathcal{I^{\prime}}$. In Figure \ref{fig:An-example-network_representation},
the path $0-1-2-4-S$ indicates that a drone starts with delivering
a package to location 1 and returns to the depot, then visits location
2 and returns to the depot again, before making its last delivery
to location 4. We define a variable $z_{ij}$ that represents whether
a drone visits location $j$ after returning from location $i$ or
not. The corresponding parameter to this variable, $c_{ij}^{d}$,
denotes the energy consumed by the drone to return empty to the depot
from location $i$ and then arrive at location $j$ with the package.
The variable $z_{0j}$ represents an outgoing arc from the dummy source
($0$), as well as the beginning of a sequence of delivery locations
visited by a drone. Therefore, adding a large penalty to the corresponding
energy consumption parameter, $c_{oj}$, we can minimize the number
of outgoing arcs (i.e., routes) from the source, and thus minimize
the required fleet size, which is the number of drones used to deliver
packages to all locations.

The mixed-integer programming (MIP) formulation of the drone deployment
optimization problem is presented below.

\subsubsection*{Objective Function}

The goal of our drone deployment optimization problem is to find the
optimal drone routing that minimizes both investment and operating
costs in delivering all customer orders within a specified planning
horizon. The investment cost includes the drone and battery costs,
whereas the operating cost includes the energy consumption, drone
operator, and maintenance costs.

\begin{eqnarray}
 & \textnormal{min}\,\,C_{E}\left(\sum_{i\in\mathcal{I^{\prime}}}\sum_{j\in\mathcal{I}^{\prime}}c_{ij}^{d}z_{ij}\right)+\left(C_{d}+\frac{C_{L}}{n_{d}}+C_{d}^{M}\right)\sum_{j\in\mathcal{I}^{\prime}}z_{0j}+C_{bat}^{d}\left(\sum_{i\in\mathcal{I^{\prime}}}y_{i}\right)\label{eq:objFunction_minRoutingCost_two_index}
\end{eqnarray}

The objective function (\ref{eq:objFunction_minRoutingCost_two_index})
seeks to minimize the total energy consumption cost (component 1),
total investment and operating costs of the drones (component 2),
and investment costs of the batteries (component 3) in delivering
all customer orders maintaining the specified pickup time windows
within the planning horizon. As the costs of drones and batteries
are functions of the number of drones required and the number of times
battery needs to be replaced, respectively, an optimal solution to
this drone routing MIP model provides the minimum required fleet size
(i.e., required number of drones) and the minimum required number
of battery replacements, which is also the minimum required additional
batteries, to deliver all customer orders satisfying the specified
time windows.

\subsubsection*{Assignment and Flow Balance Constraints}

As mentioned earlier, in this drone routing problem, each customer
order is delivered once by only one drone and the drones return to
the depot after delivering a package to each location before visiting
the next location.

\begin{eqnarray}
 & \sum_{i\in\mathcal{I}^{\prime}}z_{ij}=1\qquad\forall j\in\mathcal{I}\label{eq:minRoutingCost_cnstr_1}\\
 & \sum_{j\in\mathcal{I}^{\prime}}z_{ij}=1\qquad\forall i\in\mathcal{I}\label{eq:minRoutingCost_cnstr_2}\\
 & \sum_{i\in\mathcal{I^{\prime}}}z_{ij}=\sum_{i\in\mathcal{I}^{\prime}}z_{ji}\qquad\forall j\in\mathcal{I}\label{eq:minRoutingCost_cnstr_3}\\
 & \sum_{j\in\mathcal{I}}z_{0j}=\sum_{j\in\mathcal{I}}z_{jS}\label{eq:minRoutingCost_cnstr_4}
\end{eqnarray}

Constraints (\ref{eq:minRoutingCost_cnstr_1})\textendash (\ref{eq:minRoutingCost_cnstr_2})
ensure that each customer order is delivered once by only one drone.
Constraints (\ref{eq:minRoutingCost_cnstr_3}) and (\ref{eq:minRoutingCost_cnstr_4})
are the flow balance constraints for the customer locations, and dummy
source-sink pair, respectively, ensuring that if a drone visits a
location $i$, it must return to depot from that location to visit
the next location in the delivery sequence until the drone completes
its route.

\subsubsection*{Package Weight and Delivery Range Constraints}

As mentioned in Section \ref{sec:Problem-description}, a drone cannot
be used to deliver a package to a location that is out of its delivery
range. Constraints (\ref{eq:energy_range_cnstr}) represent this delivery
range limitation. In addition, a drone cannot be used to deliver a
customer order if the package weight exceeds the maximum weight carrying
capacity, which is enforced by constraints (\ref{eq:payload_limitation_cnstr}).

\begin{eqnarray}
 & \sum_{i\in\mathcal{I}^{\prime}}\left(E_{dHj}+E_{djH}\right)z_{ij}\leq\left(ch_{d}^{0}-ch_{d}^{min}\right)\qquad\forall j\in\mathcal{I}\label{eq:energy_range_cnstr}\\
 & M_{j}\sum_{i\in\mathcal{I}^{\prime}}z_{ij}\leq M_{d}^{max}\qquad\forall j\in\mathcal{I}\label{eq:payload_limitation_cnstr}
\end{eqnarray}

\subsubsection*{Pickup Time Window Constraints}

Constraints (\ref{eq:minRoutingCost_cnstr_5}) ensure that the pickup
time of the packages for the delivery locations are computed consistently
accounting for the sequence in which each drone visits the locations.
Constraints (\ref{eq:minRoutingCost_cnstr_6}) enforce the time windows
within which each package must be picked-up for the delivery locations
by the drones.

\begin{eqnarray}
 & f_{j}\geq f_{i}+\left(t_{\ell}+t_{Hi}+t_{u}+t_{iH}\right)+t_{d}^{bat}y_{i}-M(1-z_{ij})\qquad\forall i,j\in\mathcal{I^{\prime}},i\neq j\label{eq:minRoutingCost_cnstr_5}\\
 & e_{i}\leq f_{i}\leq\ell_{i}\qquad\forall i\in\mathcal{I}\label{eq:minRoutingCost_cnstr_6}
\end{eqnarray}

Here, $\ell_{i}=e_{i}+\Delta$ is the maximum permissible delayed
pickup time for delivery location $i$ and $M$ is a large positive
number.

\subsubsection*{Battery Replacement Constraints}

To ensure safe operation of the drones, before visiting the next location,
we need to compute how much energy would be remaining in the drone's
battery after returning to the depot from its delivery location and
replace the battery if needed. Constraints (\ref{eq:minRoutingCost_cnstr_7})
compute the potential remaining energy in the drones' batteries after
returning from the first delivery location in each delivery sequence
(i.e., route). Each constraint (\ref{eq:minRoutingCost_cnstr_8})
computes the potential remaining energy in the drone's battery after
returning from each delivery location except for the first location
in each delivery sequence. This constraint ensures the potential remaining
energy in the drone's battery is computed consistently accounting
for the sequence in which the delivery locations are visited by the
drone. For example, if location $4$ is visited after location $3$
(i.e., $z_{34}=1$), then the potential remaining energy after returning
from delivery location $4$, $g_{4}^{\prime}$, is computed by subtracting
the total energy required to arrive and return from location $4$
($E_{dH4}+E_{d4H}$) from the remaining energy after returning from
location $3$ ($g_{3}$), i.e., $g_{4}^{\prime}\leq g_{3}-\left(E_{dH4}+E_{d4H}\right).$

\begin{eqnarray}
 & g_{j}^{\prime}\leq ch_{d}^{0}-\left(E_{dHj}+E_{djH}\right)+M(1-z_{0j})\qquad\forall j\in\mathcal{I^{\prime}}\label{eq:minRoutingCost_cnstr_7}\\
 & g_{j}^{\prime}\leq g_{i}-\left(E_{dHj}+E_{djH}\right)+M(1-z_{ij})\qquad\forall i,j\in\mathcal{I^{\prime}},i\neq j\label{eq:minRoutingCost_cnstr_8}
\end{eqnarray}

If the potential remaining energy after returning from a delivery
location $j$ is less than the minimum required energy (i.e., $g_{j}^{\prime}<ch_{d}^{min}$),
then we need to replace the current battery with a fully-charged one
(i.e., $y_{i}=1$) before the drone leaves the depot to visit location
$j$. This condition is represented by constraints (\ref{eq:minRoutingCost_cnstr_9}).
But, if the potential remaining energy after returning from delivery
location $j$ is higher than the minimum required energy (i.e., $g_{j}^{\prime}\geq ch_{d}^{min}$),
then the battery must not be replaced before visiting $j$, and thus,
the drone should continue its operation with the remaining energy
($g_{i}$) after returning from its previous location $i$. Constraints
(\ref{eq:minRoutingCost_cnstr_10}) enforce this condition. Here,
$M_{u}$ and $M_{\ell}$ are the upper and lower bounds of $\left(g_{j}^{\prime}-ch_{d}^{min}\right)$,
respectively.

\begin{eqnarray}
 & g_{j}^{\prime}-ch_{d}^{min}\geq M_{\ell}y_{i}+M(z_{ij}-1)\quad\forall i\in\mathcal{I^{\prime}},j\in\mathcal{I^{\prime}\setminus\mathit{\{0\}}},i\neq j\label{eq:minRoutingCost_cnstr_9}\\
 & g_{j}^{\prime}-ch_{d}^{min}\leq M_{u}(1-y_{i})+M(1-z_{ij})\quad\forall i\in\mathcal{I^{\prime}},j\in\mathcal{I^{\prime}\setminus\mathit{\{0\}}},i\neq j\label{eq:minRoutingCost_cnstr_10}
\end{eqnarray}

Depending on whether the battery is replaced or not before visiting
the next location $j$, we need to compute the actual energy remaining
after returning from location $j$. If the current battery is replaced
with a fully-charged one after returning from location $i$ (or before
visiting location $j$), then the battery energy is resumed to the
initial energy level $ch_{d}^{0}$. Therefore, if the battery is replaced
before visiting location $j$ (i.e., $y_{i}=1$), then the energy
remaining after returning from location $j$ is computed by subtracting
the total energy consumption due to round-trip travel to location
$j$ from the initial battery energy, i.e., $g_{j}\leq ch_{d}^{0}-\left(E_{dHj}+E_{djH}\right)$.
Constraints (\ref{eq:minRoutingCost_cnstr_11}) represent this condition.
If the drone's battery is not replaced before visiting location $j$
(i.e., $y_{i}=0$), the energy remaining after returning from location
$j$ is computed based on the remaining energy ($g_{i}$) after returning
from its previous location $i$ as $g_{j}\leq g_{i}-\left(E_{dHj}+E_{djH}\right)$.
This condition is represented by constraints (\ref{eq:minRoutingCost_cnstr_12}).
Each constraint (\ref{eq:minRoutingCost_cnstr_13}) computes the actual
energy remaining after returning from the first location in each delivery
sequence. Here, $Q_{u}$ is the upper of $\left(ch_{d}^{0}-g_{j}\right)$.

\begin{eqnarray}
 & g_{j}\leq ch_{d}^{0}-\left(E_{dHj}+E_{djH}\right)+Q_{u}(1-y_{i})+M(1-z_{ij})\quad\forall i,j\in\mathcal{I^{\prime}},i\neq j\label{eq:minRoutingCost_cnstr_11}\\
 & g_{j}\leq g_{i}-\left(E_{dHj}+E_{djH}\right)+Q_{u}y_{i}+M(1-z_{ij})\quad\forall i\in\mathcal{I}^{\prime}\setminus\mathit{\{0\}},j\in\mathcal{I^{\prime}},i\neq j\label{eq:minRoutingCost_cnstr_12}\\
 & g_{j}\leq ch_{d}^{0}-\left(E_{dHj}+E_{djH}\right)+Q_{u}y_{0}+M(1-z_{0j})\quad\forall j\in\mathcal{I^{\prime}}\label{eq:minRoutingCost_cnstr_13}
\end{eqnarray}

\subsubsection*{Sign Restriction Constraints}

Constraints (\ref{eq:minRoutingCost_cnstr_14})\textendash (\ref{eq:minRoutingCost_cnstr_16})
represent the binary nature and sign restrictions of the variables.

\begin{eqnarray}
 & z_{ij}\in\{0,1\}\qquad\forall i,j\in\mathcal{I^{\prime}},i\neq j\label{eq:minRoutingCost_cnstr_14}\\
 & g_{i}\geq0\qquad\forall i\in\mathcal{I^{\prime}}\label{eq:minRoutingCost_cnstr_15}\\
 & g_{i}^{\prime}\geq-ch_{d}^{min}\qquad\forall i\in\mathcal{I^{\prime}}\\
 & f_{i}\geq0\qquad\forall i\in\mathcal{I^{\prime}}\label{eq:minRoutingCost_cnstr_16}
\end{eqnarray}

\subsection{\label{subsec:Modeling_approach_mixed_fleet}Model for Mixed Fleet
of Drones}

We propose a two-phase modeling approach to solve the drone deployment
optimization problem for a mixed fleet of drones, which contains different
types of drones in the fleet. In the first phase, we develop an integer
programming model to find the energy-efficient assignment of the different
types of drones to deliver packages to the delivery locations. This
drone assignment model (\ref{drone_assignment_model}) accounts for
the limitations of battery energy and package weight carrying capacity
for different types of drones and provides the optimal assignment
of drone types to the delivery locations.

\begin{subequations}

\begin{eqnarray}
 & \textnormal{min} & \sum_{i\in\mathcal{I}}\sum_{d\in\mathcal{D}}\left(E_{dHi}+E_{diH}\right)x_{di}\label{eq:objFun_droneAssignment}\\
 & \textnormal{s.t.,} & \sum_{d\in\mathcal{D}}x_{di}=1\qquad\forall i\in\mathcal{I}\label{eq:Cnstr_1_droneAssignment}\\
 &  & \left(E_{dHi}+E_{diH}\right)x_{di}\leq\left(ch_{d}^{0}-ch_{d}^{min}\right)\qquad\forall d\in\mathcal{D,}i\in\mathcal{I}\label{eq:Cnstr_2_droneAssignment}\\
 &  & M_{i}x_{di}\leq M_{d}^{max}\qquad\forall i\in\mathcal{I},d\in\mathcal{D}\label{eq:Cnstr_3_droneAssignment}\\
 &  & x_{di}\in\{0,1\}\qquad\forall i\in\mathcal{I},d\in\mathcal{D}\label{eq:Cnstr_4_droneAssignment}
\end{eqnarray}

\label{drone_assignment_model}

\end{subequations}

The objective function (\ref{eq:objFun_droneAssignment}) seeks to
optimally assign drone types to the delivery locations that minimizes
the total energy consumption in delivering packages to all locations.
Each constraint (\ref{eq:Cnstr_1_droneAssignment}) ensures that only
one type of drone can be assigned to deliver a package to a location.
A drone type $d$ cannot be assigned to deliver a package to a location
$i$ if that drone type consumes more energy than the maximum allowable
energy consumption (i.e., $ch_{d}^{0}-ch_{d}^{min}$) in delivering
the package to that location. Constraints (\ref{eq:Cnstr_2_droneAssignment})
enforce this battery energy limitations for the drone types. Each
constraint (\ref{eq:Cnstr_3_droneAssignment}) represents that a drone
type $d$ cannot be assigned to deliver a package to a location $i$
if the package weight for that delivery, $M_{i}$, exceeds the maximum
weight carrying capacity, $M_{d}^{max}$, of that drone type. Constraints
(\ref{eq:Cnstr_4_droneAssignment}) represent the binary nature of
the assignment decision variables.

The output of this first-phase modeling is the subsets of delivery
locations and their corresponding drone types. We use the MIP model
for homogeneous fleet of drones presented in Section (\ref{subsec:Model_identical_drones})
for each subset of delivery locations and the corresponding drone
type $d$ to obtain the required number of drones of type $d$, the
total energy consumption in delivering packages to the locations in
that subset, and the required number of battery replacements. Combining
the results from each subset and drone type, we obtain the minimum
total number of drones of different types required in the fleet, which
is the required fleet size, and the total energy consumption by all
drones in delivering all packages. A flowchart of this two-phase modeling
approach is demonstrated in Figure \ref{fig:Two-phase-modeling-approach}
(see Appendix \ref{sec:Flowchart_mixed_fleet}) to clarify the modeling
concept.

\subsection{\label{subsec:Acceleration-Techniques}Acceleration Techniques}

We observed that the number of $z_{ij}$ variables in the developed
drone deployment MIP model increases by $|\mathcal{I}^{\prime}|$
for each additional delivery location in the data. A large number
of constraints are also added in the model for each additional delivery
location. Therefore, as the number of delivery locations in our drone
deployment problem increases, the computational complexity of the
MIP model increases exponentially. To alleviate the computational
burden in solving the MIP model, we implement two problem-specific
valid inequalities that substantially reduce the number of $z_{ij}$
variables. We observed that in the optimal solution, most of the $z_{ij}$
variables are zero. Therefore, we implement the following two valid
inequalities to remove those unnecessary $z_{ij}$ variables from
the candidate solutions:
\begin{enumerate}
\item A drone cannot pick up the package for delivery location $j$ immediately
after serving location $i$, if the maximum permissible delayed pickup
time for $j$ ($\ell_{j}$) is earlier than the earliest time the
drone can return to the depot after serving location $i$. Each constraint
(\ref{eq:valid_inequality_1}) removes the $z_{ij}$ variables from
the candidate solutions for which this condition holds.
\begin{eqnarray}
 & \mathnormal{\textnormal{if\,\,}}\,e_{i}+t_{\ell}+t_{Hi}+t_{u}+t_{iH}>\ell_{j}\nonumber \\
 & \qquad\textnormal{then\,\,\,\,}z_{ij}=0\quad\forall i,j\in\mathcal{I}\label{eq:valid_inequality_1}
\end{eqnarray}
\item If the time difference between the maximum delayed pickup time for
delivery location $j$ and returning from location $i$ is longer
than the time required for a drone to serve the furthest location,
then the drone will not wait for such a long duration to pick up the
package for location $j$ after returning from location $i$. This
is because, the objective of our model is to minimize the total number
of drones required to deliver packages to all customer orders. Therefore,
the drone will deliver a package to another location $k$ after returning
from location $i$, meaning $z_{ij}=0$. Constraints (\ref{eq:valid_inequality_2})
remove the $z_{ij}$ variables from the candidate solutions for which
this condition holds.
\end{enumerate}
\begin{eqnarray}
 & \mathnormal{\textnormal{if}}\,\,\,\ell_{j}>\ell_{i}+t_{Hi}+t_{iH}+2(t_{\ell}+t_{u}+2t_{Range})\nonumber \\
 & \qquad\textnormal{then\,\,\,}z_{ij}=0\quad\forall i,j\in\mathcal{I}\label{eq:valid_inequality_2}
\end{eqnarray}

\section{Case Study Results and Managerial Insights}

We used our mathematical models described in Section \ref{sec:Model-Formulations}
to provide insights into the following questions:
\begin{enumerate}
\item How do the drone operating parameters\textemdash speed, package weight,
flight path, minimum required energy in the battery\textemdash affect
the delivery range?
\item How does the drone speed, package weight, delivery method, flight
path affect the total and average energy consumption in delivering
packages to customer locations?
\item How does the required fleet size (i.e., number of drones required)
changes with the drone speed, package weight, and flight path?
\item How sensitive is the required fleet size to changes in the pickup
time window?
\item How frequently do batteries needs to be replaced as the drone speed,
package weight, flight path, and minimum required energy change?
\item How does the use of a mixture of drones in the fleet affect the delivery
range, total energy consumption, and the required fleet size?
\end{enumerate}

\subsection{Experimental Setup}

We implemented the drone deployment optimization models described
in Section \ref{sec:Model-Formulations} using Python 3.7 with Cplex
optimizer 12.10 \citep{CPLEXOPT}. We conducted numerical experiments
using a one-hour Spatio-temporal prepared food delivery data from
the San Francisco Bay Area in California. This delivery data contains
the business provider and 126 customers\textquoteright{} geographic
locations (i.e., latitude and longitude) and the timing of when a
package is ready for pick up. We synthetically generated the package
weight for each customer order or delivery location. We conducted
drone flight tests using two rotary drones\textemdash DJI Matrice
600 Pro and Tarot 650\textemdash for different operating conditions\textemdash drone
speed, flight path (e.g., straight, over the road networks), and package
weight. DJI Matrice 600 Pro is a rotary hexacopter that has higher
package weight carrying capacity at a longer distance than the Tarot
650, which is a rotary quadcopter. Table \ref{tab:Fixed-parameter-values}
lists the parameters and their values that are kept fixed in all numerical
experiments. While designing and conducting the drone flight tests,
we used United States Customery System (USCS) units. Therefore, in
Tables \ref{tab:Fixed-parameter-values} and \ref{tab:Parameters-for-sensitivity},
we provided the parameter values both in USCS and SI units. We used
the amortized hourly costs of the drones and batteries to run our
numerical experiments.

\begin{table}[h]
\caption{\label{tab:Fixed-parameter-values}Fixed parameter values.}

\centering{}{\small{}}%
\begin{tabular}{>{\raggedright}m{4cm}>{\centering}p{2cm}>{\centering}p{2cm}>{\centering}p{2cm}>{\centering}p{2cm}}
\hline 
\multirow{2}{4cm}{\centering{}\textbf{Parameters}} & \multicolumn{2}{c}{\textbf{DJI Matrice 600 Pro Values}} & \multicolumn{2}{c}{\textbf{Tarot 650 Values}}\tabularnewline
\cline{2-5} \cline{3-5} \cline{4-5} \cline{5-5} 
 & \textbf{USCS Unit} & \textbf{SI Unit} & \textbf{USCS Unit} & \textbf{SI Unit}\tabularnewline
\hline 
\hline 
\centering{}Initial battery energy (fully-charged), $ch^{0}$  & \centering{}600 watt-hour  & \centering{}2.16 Mega-Joule (MJ) & \centering{}177.6 watt-hour & 0.639 Mega-Joule (MJ)\tabularnewline
\centering{}Package loading time, $t_{\ell}$  & \centering{}5 minutes  & \centering{}5 minutes  & \centering{}5 minutes  & \centering{}5 minutes \tabularnewline
\centering{}Package unloading time (landing), $t_{u}$  & \centering{}30 seconds  & \centering{}30 seconds  & \centering{}30 seconds  & \centering{}30 seconds \tabularnewline
\centering{}Battery replacement time, $t_{d}^{bat}$  & \centering{}5 minutes  & \centering{}5 minutes  & \centering{}5 minutes  & \centering{}5 minutes \tabularnewline
\centering{}Flying height  & \centering{}200 feet  & \centering{}60.96 meters  & \centering{}200 feet  & \centering{}60.96 meters \tabularnewline
\centering{}Hovering duration while delivering package (landing)  & \centering{}5 seconds  & \centering{}5 seconds  & \centering{}5 seconds  & \centering{}5 seconds \tabularnewline
\centering{}Hovering duration while returning to depot & \centering{}5 seconds  & \centering{}5 seconds  & \centering{}5 seconds  & \centering{}5 seconds \tabularnewline
\centering{}Cost of a drone  & \centering{}\$8000  & \centering{}\$8000  & \centering{}\$4000 & \centering{}\$4000\tabularnewline
\centering{}Cost of a drone battery  & \centering{}\$1200  & \centering{}\$1200  & \centering{}\$200 & \centering{}\$200\tabularnewline
\centering{}Labor (drone operator) cost  & \centering{}\$60/hour  & \centering{}\$60/hour  & \centering{}\$60/hour & \centering{}\$60/hour\tabularnewline
\centering{}Number of drones operated by an operator  & 5 & 5 & 5 & 5\tabularnewline
\hline 
\end{tabular}{\small\par}
\end{table}

We conducted numerical experiments for different combinations of the
parameters shown in Table \ref{tab:Parameters-for-sensitivity} to
analyze the effects these parameters on total and average energy consumption,
the required number of drones, and the required number of battery
replacements. It is to be noted that we used the same package weight
for all customer orders in each model run to demonstrate the individual
effect of package weight on the delivery range, energy consumption,
required number of drones, and battery replacements.

\begin{table}[h]
\caption{\label{tab:Parameters-for-sensitivity}Parameters for sensitivity
analysis.}

\centering{}{\small{}}%
\begin{tabular}{>{\raggedright}m{6cm}>{\centering}p{4cm}>{\centering}p{4cm}}
\hline 
\multirow{2}{6cm}{\centering{}\textbf{Parameters}} & \multicolumn{2}{c}{\textbf{Values}}\tabularnewline
\cline{2-3} \cline{3-3} 
 & \textbf{USCS Unit} & \textbf{SI Unit}\tabularnewline
\hline 
\hline 
\centering{}Drone speed & \centering{}15, 30 (mph) & \centering{}6.71, 13.41 (m/s) \tabularnewline
\centering{}Package delivery method  & \centering{}Landing, Package dropping  & \centering{}Landing, Package dropping \tabularnewline
\begin{centering}
Minimum required battery energy, $ch^{min}$ 
\par\end{centering}
(\% of initial battery energy) & \centering{}10, 15, 20, 25 & \centering{}10, 15, 20, 25 \tabularnewline
\centering{}Package weight  & \centering{}2.5, 5.0, 10.0 (lb) & \centering{}1.13, 2.27, 4.54 (Kg)\tabularnewline
\centering{}Hovering duration in package dropping  & \centering{}30, 60, 90, 120 (seconds)  & \centering{}30, 60, 90, 120 (seconds) \tabularnewline
\centering{}Pick up time window & \centering{}1, 5, 10, 15, 20, 25, 30, 35 (minutes)  & \centering{}1, 5, 10, 15, 20, 25, 30, 35 (minutes) \tabularnewline
\hline 
\end{tabular}{\small\par}
\end{table}

\begin{figure}[h]
\begin{centering}
\subfloat[\label{fig:flight_segment_delivery}Drone flight segments when delivering
the package to a customer location]{\begin{centering}
\includegraphics[width=6cm,height=6cm,keepaspectratio]{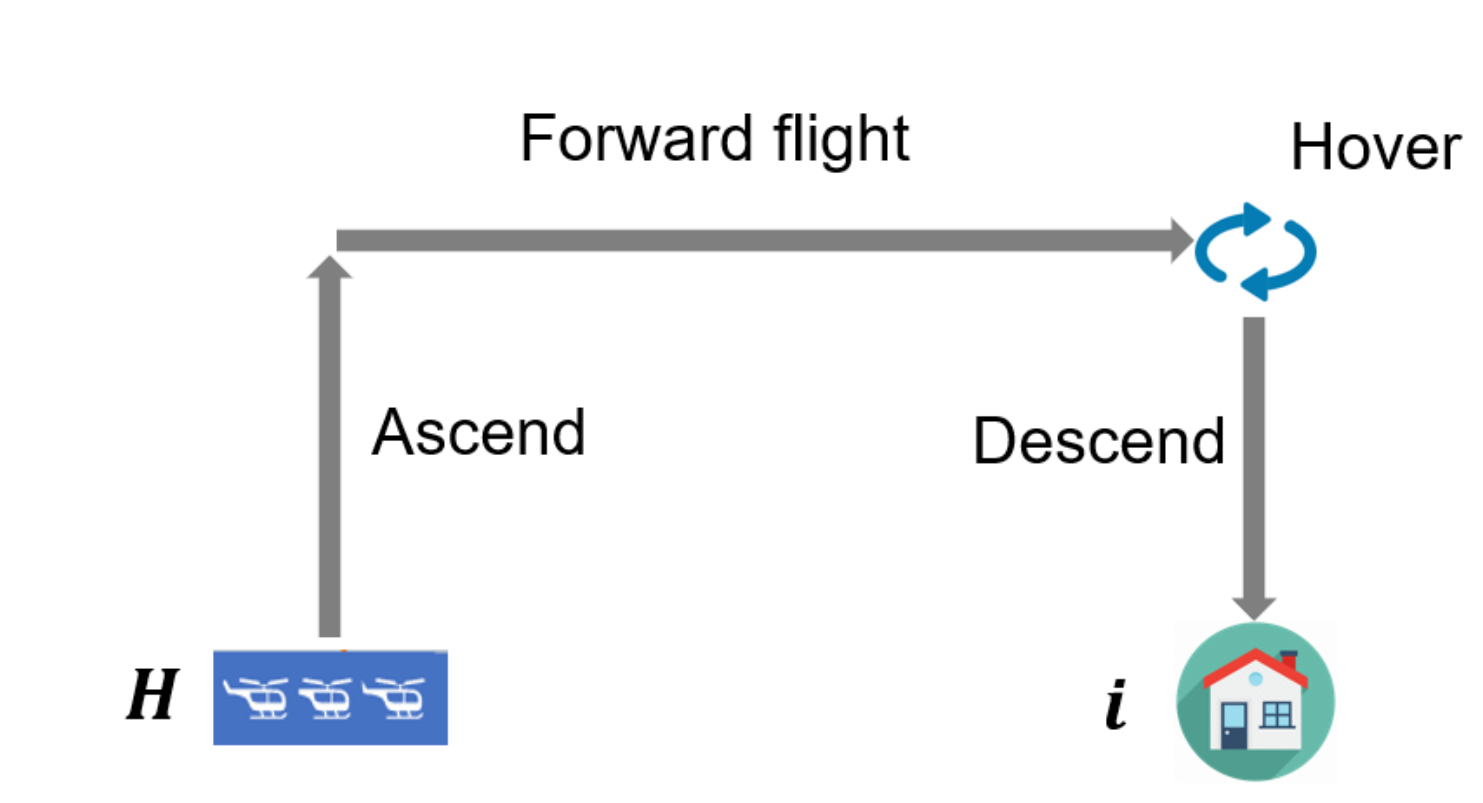}
\par\end{centering}

}~~\subfloat[\label{fig:flight-segments_return}Drone flight segments when returning
to the depot]{\begin{centering}
\includegraphics[width=6cm,height=6cm,keepaspectratio]{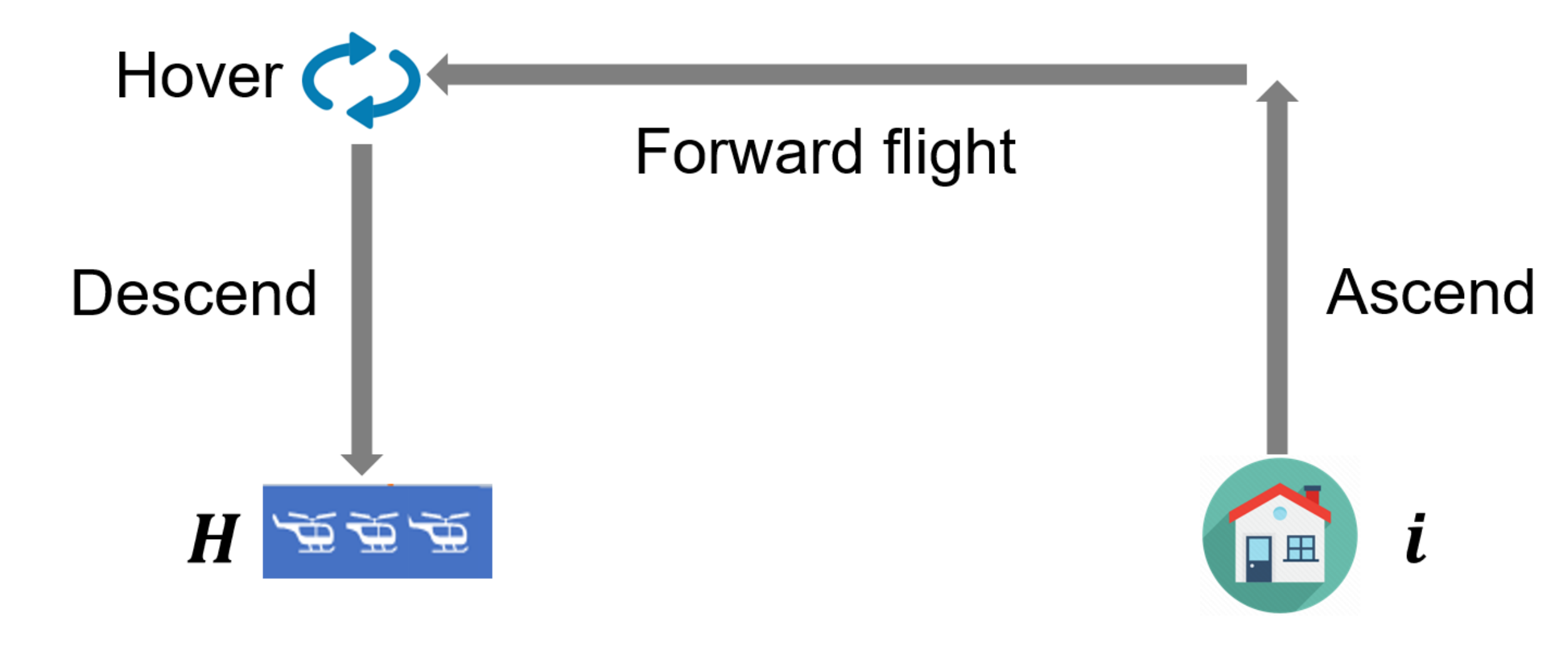}
\par\end{centering}
}
\par\end{centering}
\caption{Drone flight segments.\label{fig:Drone-flight-segments}}

\end{figure}

In the numerical experiments, we used the power consumption and flight
time data from our drone flight tests. The power consumption and flight
time data of the DJI Matrice 600 Pro drone corresponding to the flight
segments\textemdash ascend, forward flight, hover, and descend\textemdash for
different package weights and drone speeds are provided in Table \ref{tab:Energy-consumption_segments-DJI}
and \ref{tab:Time-duration_segments-DJI}, respectively in Appendix
\ref{sec:Energy-Consumption_time_segments}. Tables \ref{tab:Energy-consumption_segments-TAROT}
and \ref{tab:Time-duration_segments-TAROT} in Appendix \ref{sec:Energy-Consumption_time_segments}
show the power consumption and flight time data of the Tarot 650 drone,
respectively. Using these power consumption and flight time data,
we computed the energy consumption for each flight segment to use
in the numerical experiments. To obtain the total time required and
energy consumed by a drone to arrive with a package at a customer
location $(i)$ from the depot ($H$), we used the time and energy
of the flight segments as shown in Figure \ref{fig:flight_segment_delivery}.
Taking summation of the time of these flight segments, we obtained
the model parameter $t_{Hi}$, whereas the summation of the energy
of these segments provides the parameter $E_{dHi}$. Figure \ref{fig:flight-segments_return}
demonstrates the flight segments when the drone returns to depot ($H$)
from the delivery location $(i)$. Combining the time and energy consumption
of these segments provides the parameters $t_{iH}$ and $E_{diH}$,
respectively. Figure \ref{fig:Drone-flight-segments} demonstrates
the flight segments when the drone is delivering packages by landing
to the ground. Another delivery method is package dropping, when the
drones hover at 20 feet (6.096 meters) above the ground and drops
the package via a winch.

\subsection{Effect of Drone Speed and Package Weight on Delivery Range, Energy
Consumption, Required Fleet Size, and Battery Replacements}

The delivery range of a drone is limited by the drone's battery energy
and the fact that we need to ensure a minimum required energy remains
in the drone battery upon returning to the depot for its safe operation.
We conducted numerical experiments to evaluate the effect of drone
speed and package weight on the delivery range. In this analysis,
we used the DJI Matrice 600 Pro drone for which we used the minimum
required battery energy as 324 Kilo-Joule (KJ), which is 15\% of the
initial battery energy of 2.16 Mega-Joule (MJ). Therefore, the delivery
locations for which the DJI drone consumes more energy than 1.836
MJ ($=2.16-0.324$) to deliver the package and return, are considered
to be outside the drone delivery range. As mentioned earlier, the
delivery range is defined by the number of locations that are outside
the range of drone delivery. Table \ref{tab:speed_weight_delivery_range}
demonstrates the number of locations that are outside the drone range
for two different drone speeds and three different package weights.

\begin{table}[h]

\caption{\label{tab:speed_weight_delivery_range}Number of locations outside
drone delivery range.}

\begin{centering}
{\small{}}%
\begin{tabular}{>{\raggedright}m{3cm}>{\centering}p{1.5cm}>{\centering}p{1.5cm}>{\centering}p{1.5cm}}
\hline 
\multirow{2}{3cm}{\centering{}\textbf{Drone Speed (m/s)}} & \multicolumn{3}{c}{\textbf{Package Weight (Kg)}}\tabularnewline
\cline{2-4} \cline{3-4} \cline{4-4} 
 & \centering{}\textbf{\small{}1.13} & \centering{}\textbf{\small{}2.27} & \centering{}\textbf{\small{}4.54}\tabularnewline
\hline 
\hline 
\centering{}13.41 & 0 & 0 & 0\tabularnewline
\centering{}6.71 & \centering{}27 & \centering{}29 & 45\tabularnewline
\hline 
\end{tabular}{\small\par}
\par\end{centering}
\end{table}

We see from Table \ref{tab:speed_weight_delivery_range} that the
DJI Matrice 600 Pro can deliver packages to all 126 customer locations
when flying at 13.41 meter per second (m/s) speed. However, the drone
cannot deliver packages to all locations when flying at 6.71 m/s.
This is because, despite the average power consumption is less at
a slower speed (6.71 m/s), the drone stays in the air for a much longer
duration at 6.71 m/s speed to travel the same distance, resulting
in a much larger total energy consumption as compared to 13.41 m/s
speed in delivering a package. Due to this much larger consumption
of energy, it is not possible for the drone to safely deliver packages
to far away locations and return to the depot while flying at 6.71
m/s speed. The drone energy consumption increases as the package weight
increases. Therefore, we see from Table \ref{tab:speed_weight_delivery_range}
that more locations fall outside the delivery range when the weights
of all customer orders increases from 1.13 Kg to 2.27 Kg and 4.54
Kg. As we see that 45 locations are outside the delivery range of
the DJI drones with 6.71 m/s speed, we excluded these 45 locations
for the other analysis in this subsection to ensure a fair comparison
in the required fleet size (i.e., required number of drones), required
number of battery replacements, and the total and average energy consumption
between these two speeds.

\begin{figure}[h]
\begin{centering}
\includegraphics[width=12cm,height=12cm,keepaspectratio]{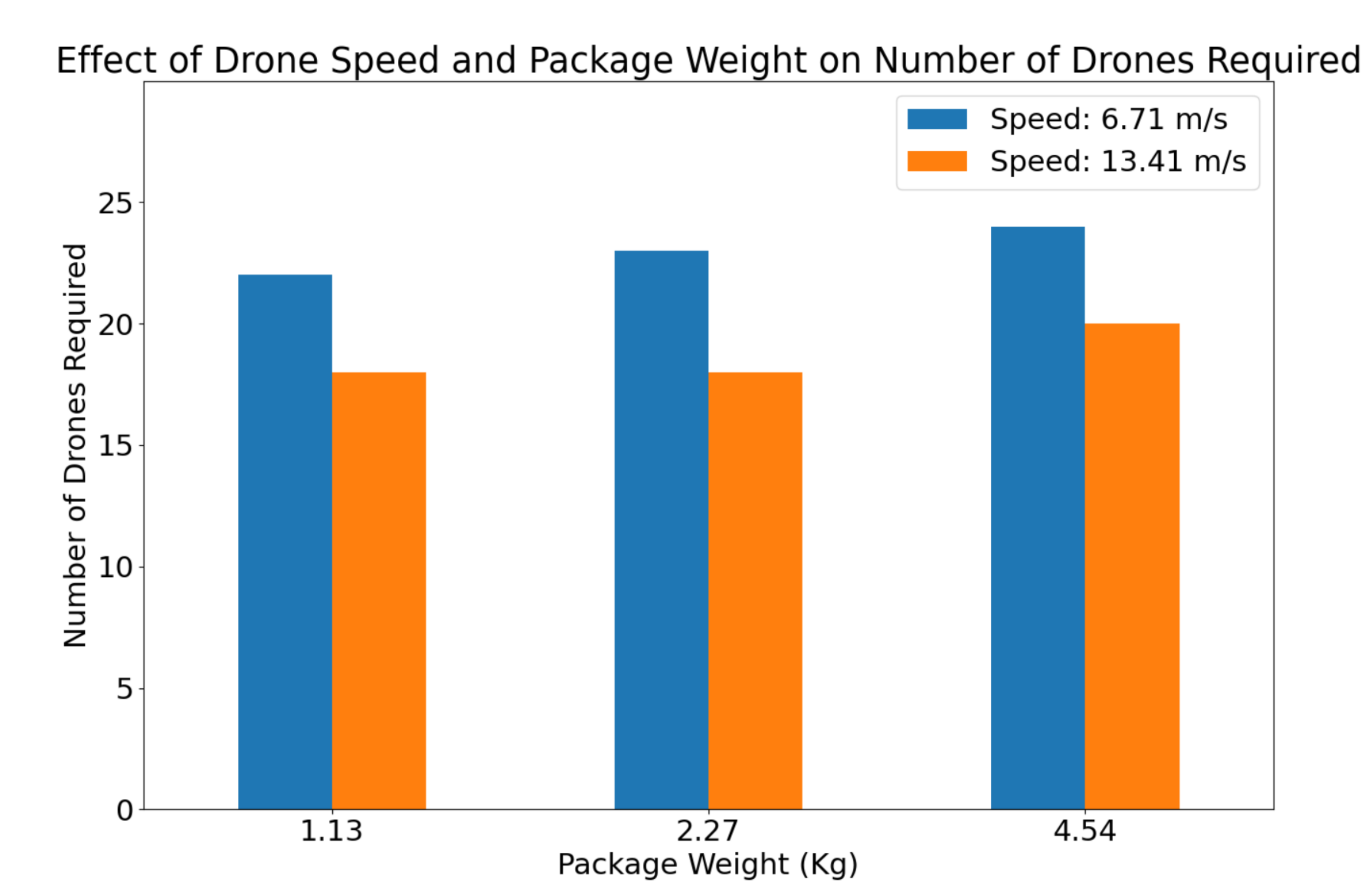}
\par\end{centering}
\caption{\label{fig:speed_weight_on_Number_of_drones}Effect of drone speed
and package weight on the required number of drones.}
\end{figure}

Figure \ref{fig:speed_weight_on_Number_of_drones} demonstrates the
variation in the required number of drones to deliver all customer
orders as the speed changes from 13.41 m/s to 6.71 m/s. In this experiments,
we used the pickup time window to be 15 minutes for the customer orders.
We see that across all three package weights, more drones are needed
with 6.71 m/s speed as compared to 13.41 m/s. The percentage increase
in the required number of drones are 22.2\%, 27.7\%, and 20\% corresponding
to package weights of 1.13 Kg, 2.27 Kg, and 4.54 Kg, respectively,
as the drone speed reduces from 13.41 m/s to 6.71 m/s. This is because,
with a slower speed (i.e., 6.71 m/s), it is not possible to maintain
the same pickup time window with the same number of drones as in 13.41
m/s. However, we see that the package weight does not have much impact
on the required number of drones. This is because, drone speed and
travel time are not affected by the package weight. One or two additional
drones are required as the package weight increases from 1.13 Kg due
to the fact that a higher number of battery replacements are needed
with a larger package weight. Due to this additional battery replacement
time, additional drones are needed to maintain the same time window.

\begin{figure}[h]

\subfloat[\label{fig:Total-energy-consumption}Total energy consumption.]{\begin{centering}
\includegraphics[width=8cm,height=8cm,keepaspectratio]{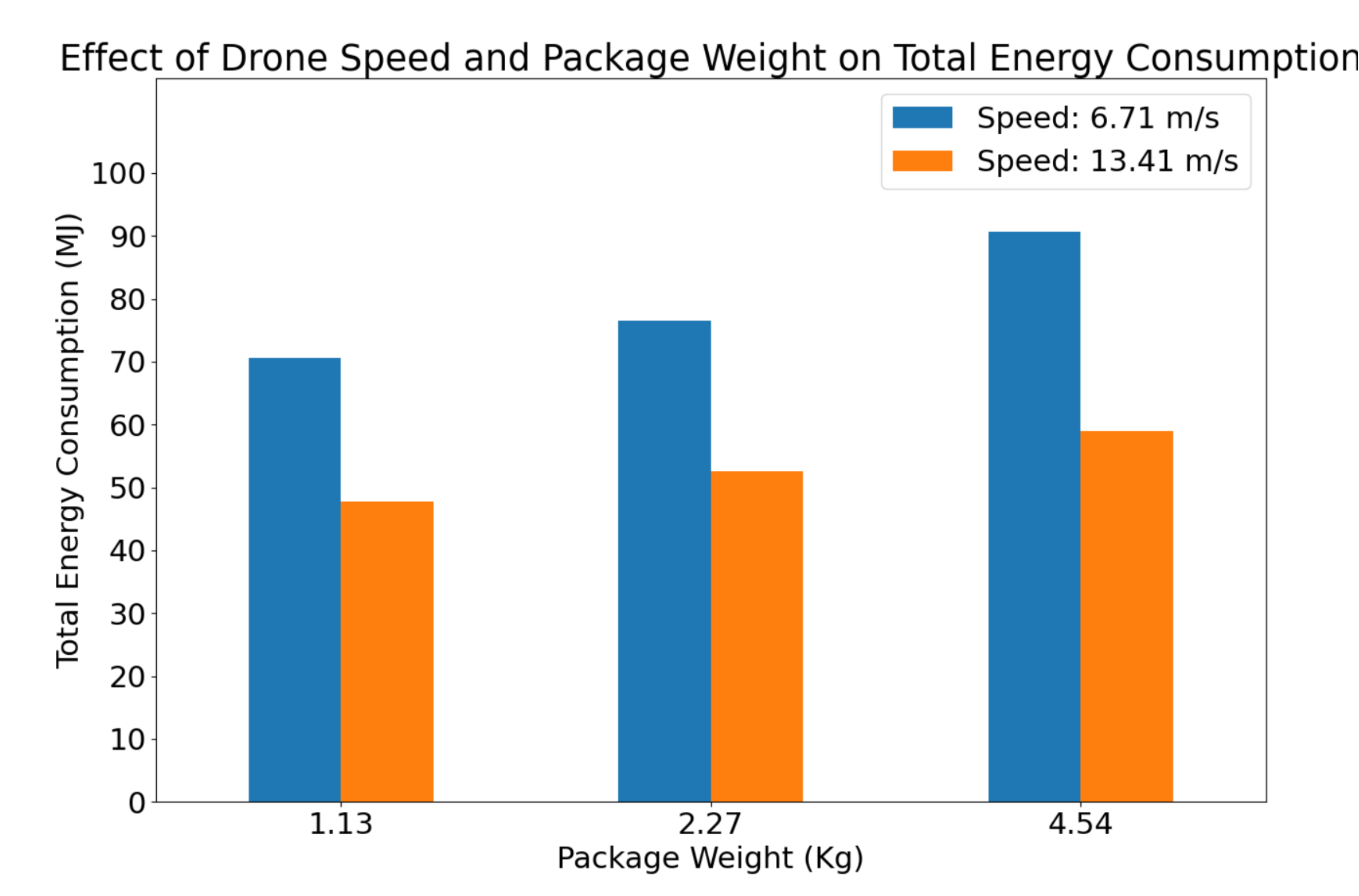}
\par\end{centering}
}~~~\subfloat[\label{fig:Average-energy-consumption}Average energy consumption.]{\begin{centering}
\includegraphics[width=8cm,height=8cm,keepaspectratio]{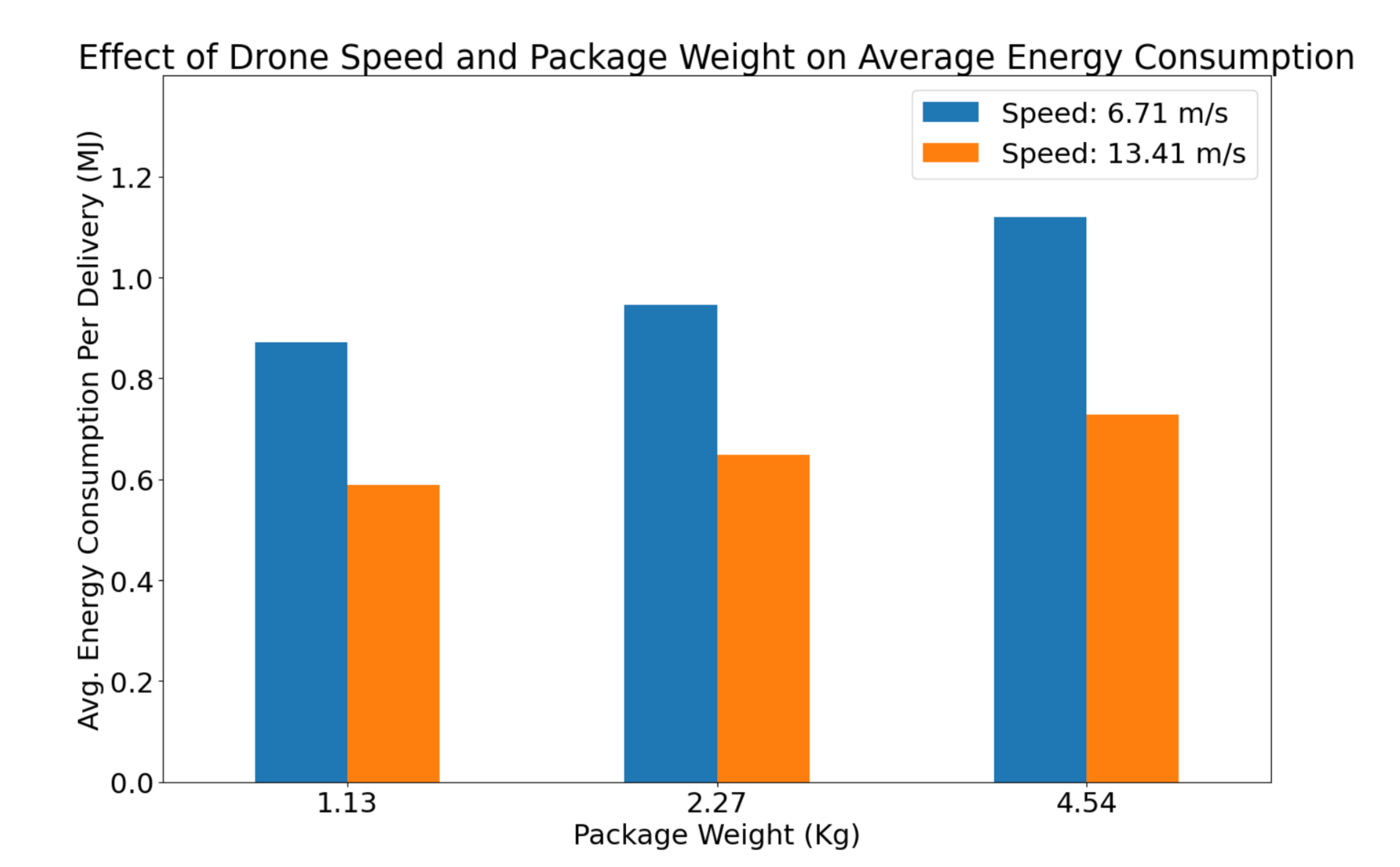}
\par\end{centering}
}\caption{\label{fig:Effect-of_speed_weight_on_energy}Effect of drone speed
and package weight on energy consumption.}

\end{figure}

Figures \ref{fig:Total-energy-consumption} and \ref{fig:Average-energy-consumption}
demonstrate the effect of drone speed and package weight on the total
energy consumption in delivering all customer orders and on the average
energy consumption per delivery, respectively. We see that the total
and average energy consumption increases as the drone speed reduces
from 13.41 m/s to 6.71 m/s across all three package weights. As mentioned
earlier, the average power consumption by a drone is lower at a slower
speed, but the drone stays in the air for a longer duration to travel
the same distance at a slower speed. Therefore, compared to the 13.41
m/s speed, flying at 6.71 m/s results in a much larger energy consumption
in delivering each customer order, eventually increasing the total
and average energy consumption. We have observed that the total energy
consumption increases by 47.96\%, 45.78\%, and 53.82\% corresponding
to the package weights of 1.13 Kg, 2.27 Kg, and 4.54 Kg, respectively,
as the drone speed reduces from 13.41 m/s to 6.71 m/s. The percentage
changes in the average energy consumption per delivery also demonstrates
the same trend and magnitude. It is evident from Figures \ref{fig:Total-energy-consumption}
and \ref{fig:Average-energy-consumption} that the energy consumption
increases as the drones carry larger package weight at both speed
levels, which is very intuitive. The percentage increase in total
energy consumption are 10.03\% and 23.45\% as the package weight of
all customer orders increases from 1.13 Kg to 2.27 Kg and 4.54 Kg,
respectively, at 13.41 m/s drone speed. We see a similar effect with
an increase in the package weight on energy consumption when the drone
speed is 6.71 m/s.

\begin{figure}[h]

\begin{centering}
\includegraphics[width=12cm,height=12cm,keepaspectratio]{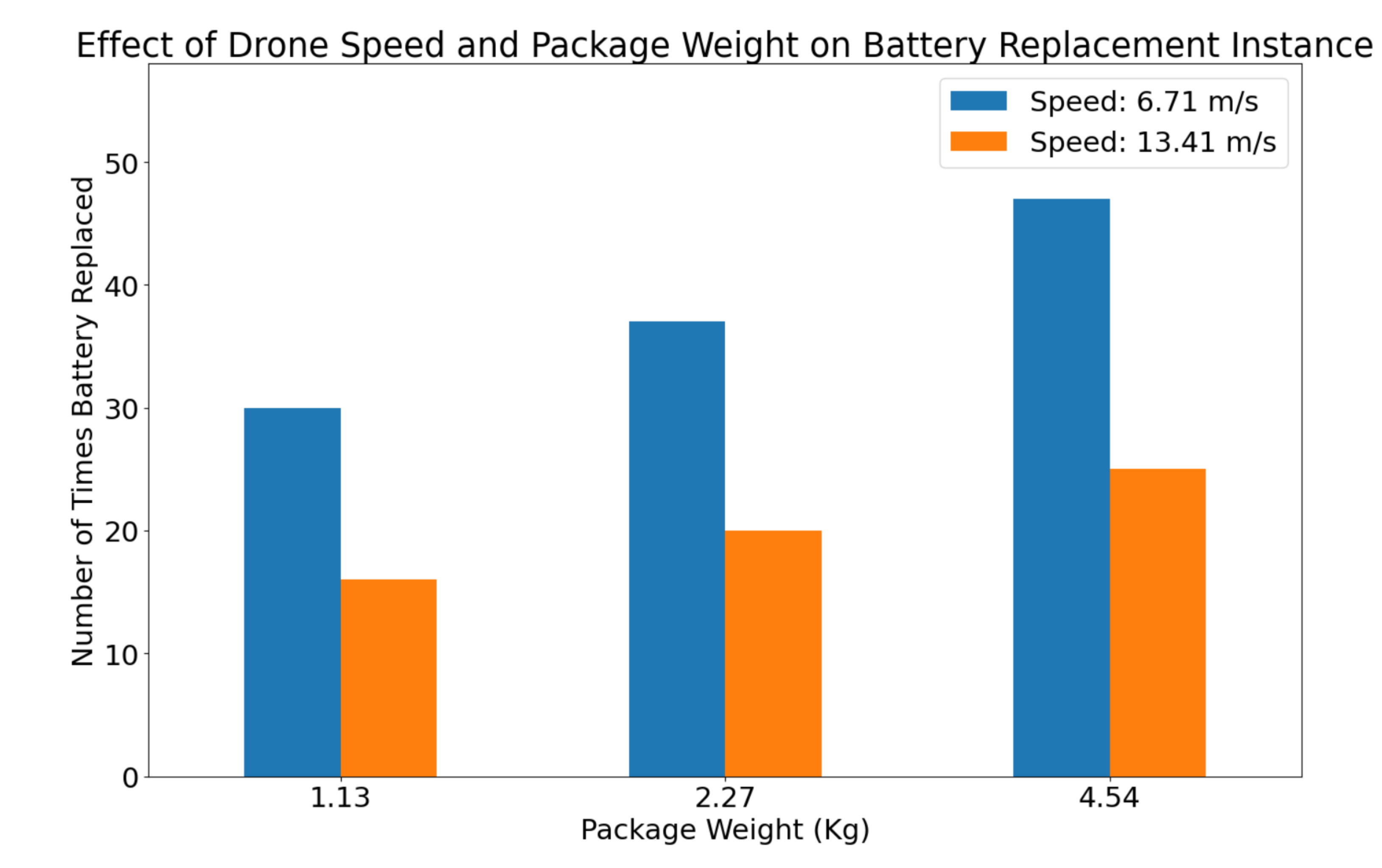}
\par\end{centering}
\caption{\label{fig:Effect_of_speed_weight_on_battery}Effect of drone speed
and package weight on the required number of battery replacements.}
\end{figure}

Figure \ref{fig:Effect_of_speed_weight_on_battery} demonstrates the
effect of drone speed and package weight on the required number of
battery replacements in delivering all 81 customer orders. We see
that the drone battery needs to be replaced a substantially larger
number of times in delivering all customer orders at 6.71 m/s speed
as compared to 13.41 m/s speed. This is due to an increased energy
consumption at 6.71 m/s speed as discussed earlier. As the energy
consumption at 6.71 m/s speed is higher than at 13.41 m/s speed in
delivering a customer order, the minimum required battery energy is
reached more frequently at a slower speed, thus necessitating more
battery replacements.

The required number of battery replacements also increases as the
package weight increases at the same speed level. Table \ref{tab:PackageWeight_on_percentage_energy_consumption}
demonstrates the percentage increase in the number of battery replacements
as the package weight of all customer orders increases from 1.13 kg
for the two speed levels.

\begin{table}[h]

\caption{\label{tab:PackageWeight_on_percentage_energy_consumption}Effect
of package weight on the number of battery replacements.}

\centering{}{\small{}}%
\begin{tabular}{>{\centering}p{3.5cm}>{\centering}p{3.5cm}>{\centering}p{3.5cm}}
\hline 
\multirow{2}{3.5cm}{\centering{}\textbf{Package Weight (Kg)}} & \multicolumn{2}{c}{\textbf{Percentage Increase in Number of Battery Replacements }}\tabularnewline
\cline{2-3} \cline{3-3} 
 & \textbf{Speed: 13.41 m/s} & \textbf{Speed: 6.71 m/s}\tabularnewline
\hline 
\hline 
\centering{}2.27 & \centering{}{\small{}25.0 } & \centering{}{\small{}23.33 }\tabularnewline
\centering{}4.54 & \centering{}{\small{}56.25 } & \centering{}{\small{}56.67 }\tabularnewline
\hline 
\end{tabular}{\small\par}
\end{table}

\subsection{Effect of Delivery Method on Energy Consumption}

While delivering a package to a customer location, the drone can either
land on the ground, as shown in Figure \ref{fig:Drone-flight-segments}
or hover 20 feet (6.096 meters) above the ground and drop the package
via a winch. As the hovering operation is energy-intensive, we compare
the energy consumption between these two delivery methods\textemdash package
dropping and landing with a package to the ground. Depending on the
type of product, the hovering duration while package dropping can
be different, ranging from 30 seconds to several minutes. Therefore,
we compared the total energy consumption in delivering all customer
orders via landing with four different package dropping methods where
the hovering durations are 30 seconds, 1 minute, 1.5 minutes, and
2 minutes. In this analysis, we used the data for the DJI Matrice
600 Pro drone flying at a speed of 13.41 m/s while carrying a package
weight of 1.13 Kg.

We would like to note that compared to landing, the drones save some
amount of ascend and descend energy while package dropping by not
descending to the ground, as well as not ascending from the ground
while returning. Table \ref{tab:Comparison_delivery_methods} demonstrates
the percentage increase in the total energy consumption for using
the different package dropping methods as compared to landing. We
see that when the hovering duration in package dropping is short (e.g.,
30 seconds), the percentage increase in the total energy consumption
is not significant compared to energy consumption in landing. As such,
the package dropping also saves some energy by avoiding to land and
ascend from the ground. However, as the hover duration becomes longer,
drones consume a much higher amount of energy while delivering a package
via dropping as compared to landing.

\begin{table}[h]
\caption{\label{tab:Comparison_delivery_methods}Comparison of energy consumption
between landing and package dropping.}

\centering{}{\small{}}%
\begin{tabular}{>{\centering}p{3.5cm}>{\centering}p{3.5cm}>{\raggedright}m{4.5cm}}
\hline 
\multirow{2}{3.5cm}{\centering{}\textbf{Delivery Method}} & \multirow{2}{3.5cm}{\centering{}\textbf{Hovering Duration}} & \multirow{2}{4.5cm}{\centering{}\textbf{\% Increase in Energy Consumption }}\tabularnewline
 &  & \tabularnewline
\hline 
\hline 
\centering{}Package dropping 1 & \centering{}30 seconds & \centering{}1.45\tabularnewline
\centering{}Package dropping 2 & 1 minute & \centering{}6.0\tabularnewline
\centering{}Package dropping 3 & 1.5 minutes & \centering{}12.48\tabularnewline
\centering{}Package dropping 4 & \centering{}2 minutes & \centering{}14.98\tabularnewline
\hline 
\end{tabular}{\small\par}
\end{table}

\subsection{Effect of Minimum Required Battery Energy on Delivery Range and Battery
Replacements}

We mentioned earlier that a minimum specified amount of energy should
be remaining in the drone battery each time a drone returns to the
depot to ensure the safe operation of the drones. This is because,
if the remaining energy in the drone battery drops below a critical
level, the drone will be grounded while flying. This critical minimum
required amount of energy for our DJI Matrice 600 Pro is 10\% of the
initial battery energy of 2.16 MJ. However, business owners may choose
a larger value of this minimum required energy than 10\% to ensure
the safe operation of the drones under adverse conditions, such as
flying in the opposite direction of strong wind. Different business
owners can have different levels of risk preferences and set the minimum
required amount of energy to different levels. Therefore, we analyzed
the effect of different values of this minimum required energy on
the delivery range and the required number of battery replacements
in delivering all customer orders.

Table \ref{tab:minRequiredEnergy_deliveryRange} demonstrates the
number of delivery locations out of drone range for four different
values of minimum required levels of battery energy\textemdash 10\%,
15\%, 20\%, and 25\% of the initial battery energy\textemdash and
three different package weights\textemdash 1.13 Kg, 2.27 Kg, and 4.54
Kg. We see that more locations fall outside the delivery range as
the minimum required battery energy increases with the same package
weight. This is because as the requirement for remaining energy in
the drone battery after returning from each delivery location increases,
this increased threshold cannot be met for a higher number of locations.
Therefore, with this increased minimum required amount of energy,
the number of customer locations that drones cannot safely deliver
packages to increases. In addition, the number of locations outside
the delivery range increases as the package weight increases within
the same minimum required battery energy level. This is because the
energy consumption increases with the package weight; thus, the minimum
required energy threshold cannot be met for a larger number of locations
and fall outside the delivery range.

\begin{table}[h]
\caption{\label{tab:minRequiredEnergy_deliveryRange}Effect of minimum required
battery energy on delivery range. $\mathbf{m}$ denotes the package
weight.}

\begin{centering}
{\small{}}%
\begin{tabular}{>{\raggedright}m{5cm}>{\centering}p{2.15cm}>{\centering}p{2.15cm}>{\centering}p{2.15cm}}
\hline 
\multirow{2}{5cm}{\centering{}\textbf{Minimum Required Battery Energy (KJ)}} & \multicolumn{3}{c}{\textbf{Number of Locations Outside Delivery Range}}\tabularnewline
\cline{2-4} \cline{3-4} \cline{4-4} 
 & \centering{}\textbf{$\mathbf{m=1.13}$ Kg} & \centering{}\textbf{$\mathbf{m=2.27}$ Kg} & \centering{}\textbf{$\mathbf{m=4.54}$ Kg}\tabularnewline
\hline 
\hline 
\centering{}216 (10\%) & 0 & 0 & 6\tabularnewline
\centering{}324 (15\%) & 0 & 0 & 10\tabularnewline
\centering{}432 (20\%) & 0 & 7 & 12\tabularnewline
\centering{}540 (25\%) & \centering{}2 & \centering{}10 & 14\tabularnewline
\hline 
\end{tabular}{\small\par}
\par\end{centering}
\end{table}

\begin{figure}[h]

\begin{centering}
\includegraphics[width=12cm,height=12cm,keepaspectratio]{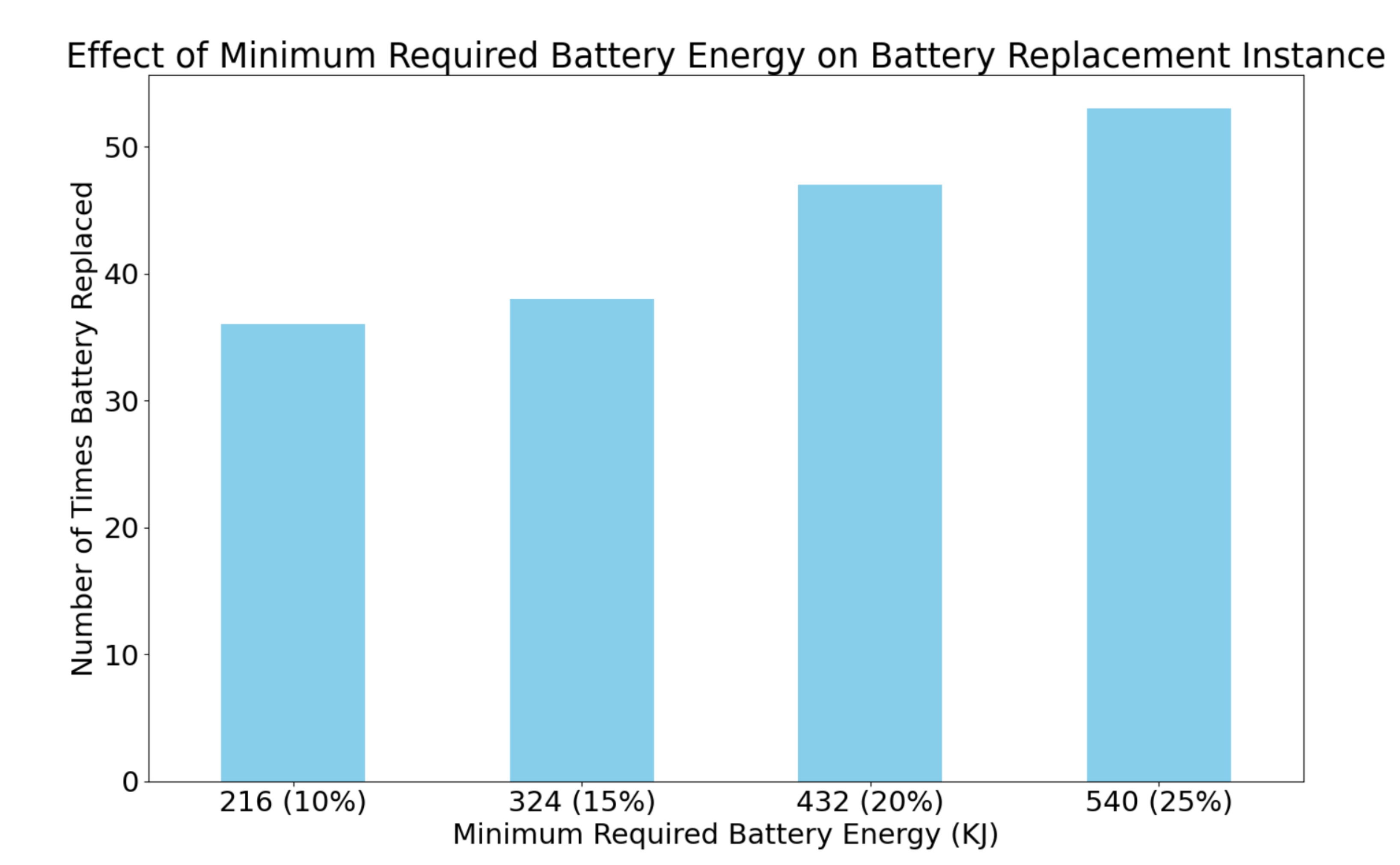}
\par\end{centering}
\caption{\label{fig:minimum_battery_energy_battery_replacements}Effect of
minimum required battery energy on the required number of battery
replacements.}
\end{figure}

This minimum required amount of battery energy also affects the required
number of battery replacements in delivering packages to customer
locations. To demonstrate this effect, we conducted numerical analysis
for the four different values of minimum required battery energy with
1.13 Kg package weight and excluding the two locations outside the
delivery range, as shown in Table \ref{tab:minRequiredEnergy_deliveryRange}.
Figure \ref{fig:minimum_battery_energy_battery_replacements} demonstrates
that the required number of battery replacements in delivering packages
to all 124 customer locations increases as the minimum required battery
energy increases. As the minimum required battery energy is set to
a higher level, the minimum threshold is reached more quickly during
drone operations, which increases the frequency of replacing the drone
batteries. We see from Figure \ref{fig:minimum_battery_energy_battery_replacements}
that the required number of battery replacements increases substantially
(i.e., 34\%) as the minimum required energy increases from 10\% to
25\% of the initial battery energy.

\subsection{Effect of Time Window on the Required Fleet Size}

We analyzed the impact of a pickup time window requirement of customer
orders on the minimum required drone fleet size (i.e., required number
of drones) to deliver all customer orders. As mentioned in Section
\ref{sec:Problem-description}, we considered a pickup time window
in which a customer wants a drone to pick the product up for delivery
to them. This time window requirement is different for different products
and business types (e.g., 2 minutes for Starbucks coffee vs. 30 minutes
for grocery items from Walmart). Our drone deployment optimization
model accounts for this time window preference of customers while
finding the optimal solution\textemdash minimum required fleet size,
required number of battery replacements, and energy consumption.

\begin{figure}[h]
\begin{centering}
\includegraphics[width=12cm,height=12cm,keepaspectratio]{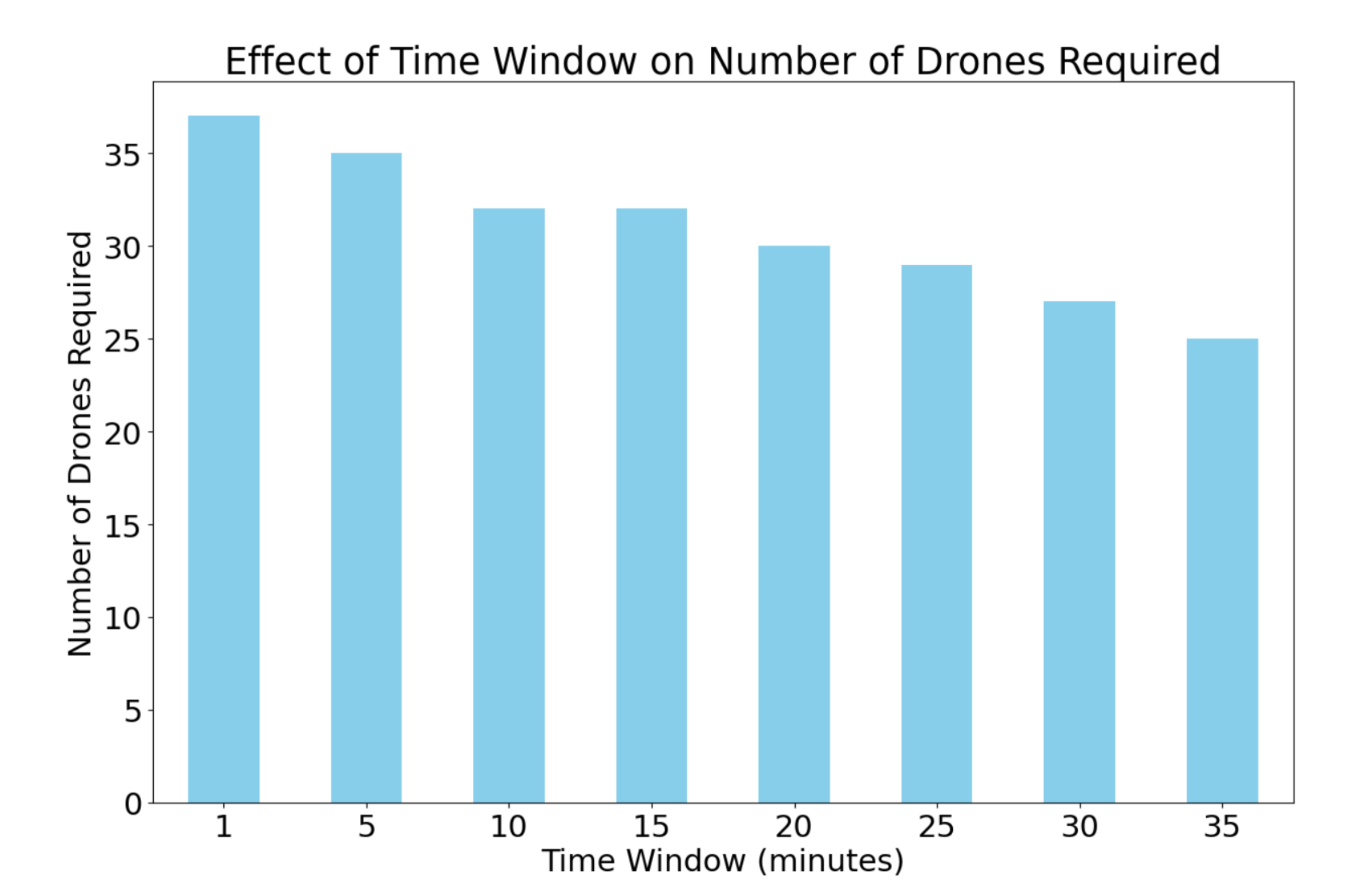}
\par\end{centering}
\caption{\label{fig:time_window_number_of_drones}Effect of time window on
the required number of drones.}

\end{figure}

Figure \ref{fig:time_window_number_of_drones} demonstrates the variation
in the required number of drones as the pickup time window increases
from 1 to 35 minutes. We see that the required number of drones in
delivering all customer orders decreases as the time window increases.
In a real-life business environment, multiple customer orders are
placed at the same time or their order placement times very close
to each other, which is reflected in the real-life delivery data we
used in our case study. In this real-life business environment, with
a more flexible pickup time window, the same drone can be used to
deliver packages to a larger number of customer orders by meeting
the time window. But, while satisfying a very short time window, it
is not possible for a drone to pick up the packages for multiple customer
orders having order placement time close to each other. Therefore,
the required number of drones to deliver all customer orders increases
as the pickup time window decreases. For instance, the required number
of drones (i.e., fleet size) increases by 48\% as the time window
reduces from 35 minutes to 1 minute.

Therefore, depending on the business and product type, the minimum
required drone fleet size can vary. Business owners can use our model
on their delivery data to determine the minimum required drone fleet
size based on the time window requirements of their business and product
types.

\subsection{Effect of Drone Flight Path on the Delivery Range, Energy Consumption,
Fleet Size, and Battery Replacements}

In real-life operations of drones, it may not always be possible to
fly the drones in a straight path from the business provider's location
to customer locations due to the presence of no-fly zones (e.g., schools,
large buildings, and restricted private property) in the flight path.
Flying over the road networks is a viable option to avoid these no-fly
zones. But, as compared to a straight path, flying over the road networks
requires the drones to travel longer distances and consume higher
amounts of energy in delivering customer orders. Our hypothesis is
that this increased distance and energy consumption affects drone
delivery range, as well as the required fleet size (i.e., number of
drones) and the number of battery replacements. Therefore, using our
drone deployment optimization model, we evaluated the impact of drone
flight path (i.e., straight vs. over the road networks) on the delivery
range, minimum number of drones and battery replacements needed.

To incorporate the road network flight path of the drones in our optimization
model, we used the Open Source Routing Machine (OSRM) API \citep{OSRM2020}
to obtain the road distance traveled and the number of turns a drone
makes to change direction while traveling to each customer location.
Based on our DJI Matrice 600 Pro drone flight test, we consider that
a drone hovers 7 seconds on average to change direction (including
the time for deceleration and acceleration) while flying over the
road networks. Using these road distances and turns, we compute the
energy consumed by the DJI Matrice 600 Pro drones and travel times
for all customer locations that are used in the optimization model.
We used the minimum required battery energy as 324 Kilo-Joule (KJ),
which is 15\% of the initial battery energy of 2.16 Mega-Joule (MJ)
in all the analyses presented in this subsection, where we used a
15 minute pickup time window.

Table \ref{tab:flight_path_on_delivery_range} shows the number of
locations out of delivery range of the DJI Matrice 600 Pro for two
different flight paths (i.e., straight and over the road networks),
two different speed levels (i.e., 6.71 m/s and 13.41 m/s), and three
different package weights (i.e., 1.13 Kg, 2.27 Kg, and 4.54 Kg). We
see that within the same package weight and speed level, a larger
number of customer locations are out of the drone delivery range as
the drones fly over the road networks as compared to the straight
path. This is because the drones consume substantially higher energy
while flying over the road networks. Therefore, given the battery
energy limitation, it is not possible for the drones to safely deliver
packages and return to the depot for a larger number of locations.
Compared to 13.41 m/s, flying over the road networks at the slower
speed (i.e., 6.71 m/s) substantially limits the drone delivery range.
For instance, 82 of the 126 locations fall outside the delivery range
when drones fly over the road networks at 6.71 m/s speed. Therefore,
in the remaining analysis in this subsection, we excluded these 82
locations from the delivery data to compare the results for both speed
levels and flight paths.

\begin{table}[h]
\caption{\label{tab:flight_path_on_delivery_range}Effect of flight path on
delivery range. $\mathbf{m}$ denotes the package weight.}

\begin{centering}
{\small{}}%
\begin{tabular}{>{\raggedright}m{3cm}>{\raggedright}m{3cm}>{\centering}p{2.15cm}>{\centering}p{2.15cm}>{\centering}p{2.15cm}}
\hline 
\multirow{2}{3cm}{\centering{}\textbf{Drone Speed (m/s)}} & \multirow{2}{3cm}{\centering{}\textbf{Flight Path}} & \multicolumn{3}{c}{\textbf{Number of Locations Outside Delivery Range}}\tabularnewline
\cline{3-5} \cline{4-5} \cline{5-5} 
 &  & \centering{}\textbf{$\mathbf{m=1.13}$ Kg} & \centering{}\textbf{$\mathbf{m=2.27}$ Kg} & \centering{}\textbf{$\mathbf{m=4.54}$ Kg}\tabularnewline
\hline 
\hline 
\multirow{2}{3cm}{\centering{}13.41} & \centering{}Straight & \centering{}0 & \centering{}0 & \centering{}10\tabularnewline
 & \centering{}Over Road Network & \centering{}15 & \centering{}23 & \centering{}37\tabularnewline
\hline 
\multirow{2}{3cm}{\centering{}6.71} & \centering{}Straight & \centering{}27 & \centering{}29 & \centering{}45\tabularnewline
 & \centering{}Over Road Network & \centering{}49 & \centering{}57 & \centering{}82\tabularnewline
\hline 
\end{tabular}{\small\par}
\par\end{centering}
\end{table}

\begin{figure}[h]
\begin{centering}
\subfloat[\label{fig:flight_path_speed_30_number_of_drones}Required number
of drones for drone speed of 13.41 m/s]{\begin{centering}
\includegraphics[width=8cm,height=8cm,keepaspectratio]{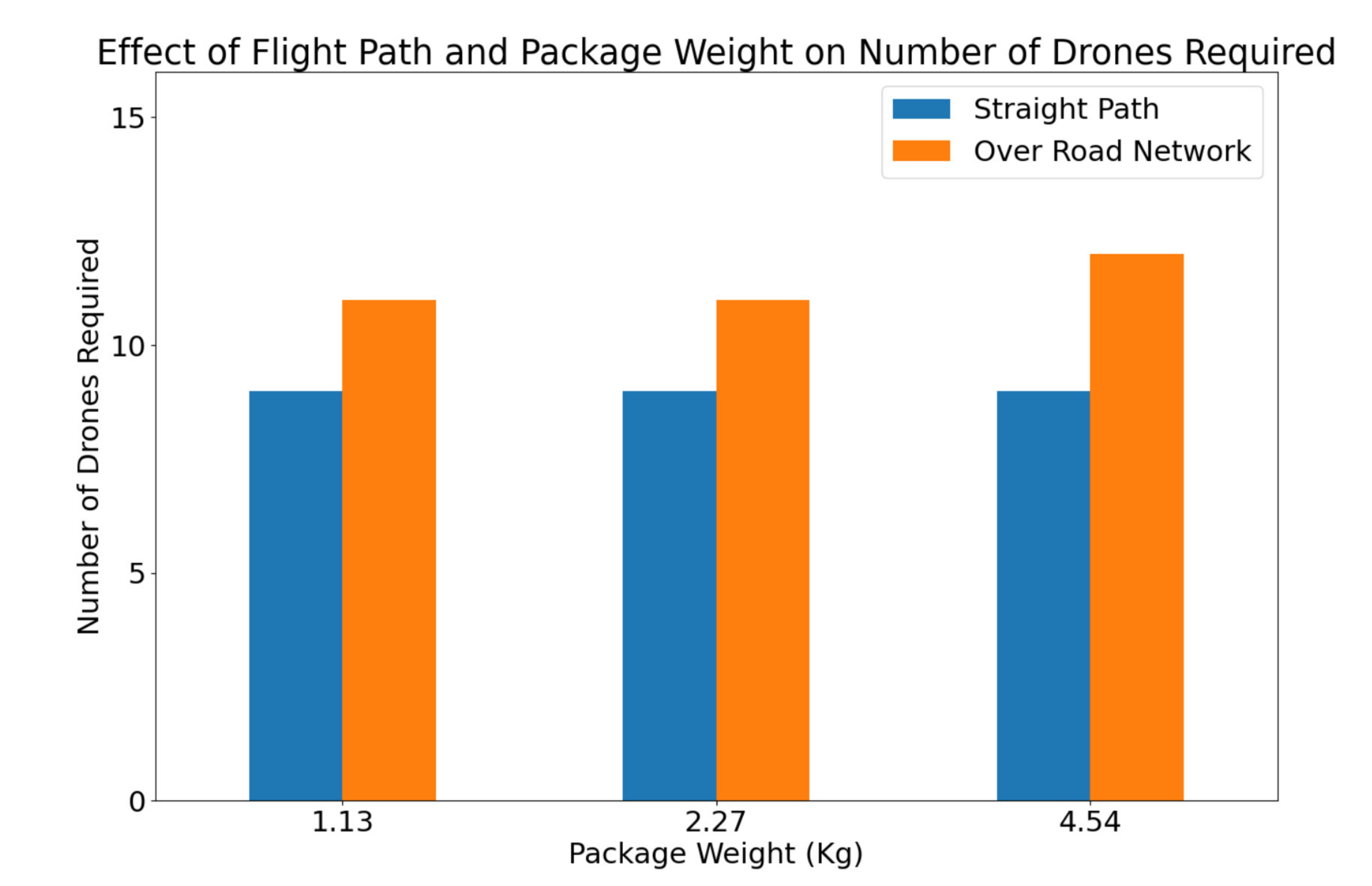}
\par\end{centering}

}~~\subfloat[\label{fig:flight_path_speed_15mph_number_of_drones}Required number
of drones for drone speed of 6.71 m/s]{\centering{}\includegraphics[width=8cm,height=8cm,keepaspectratio]{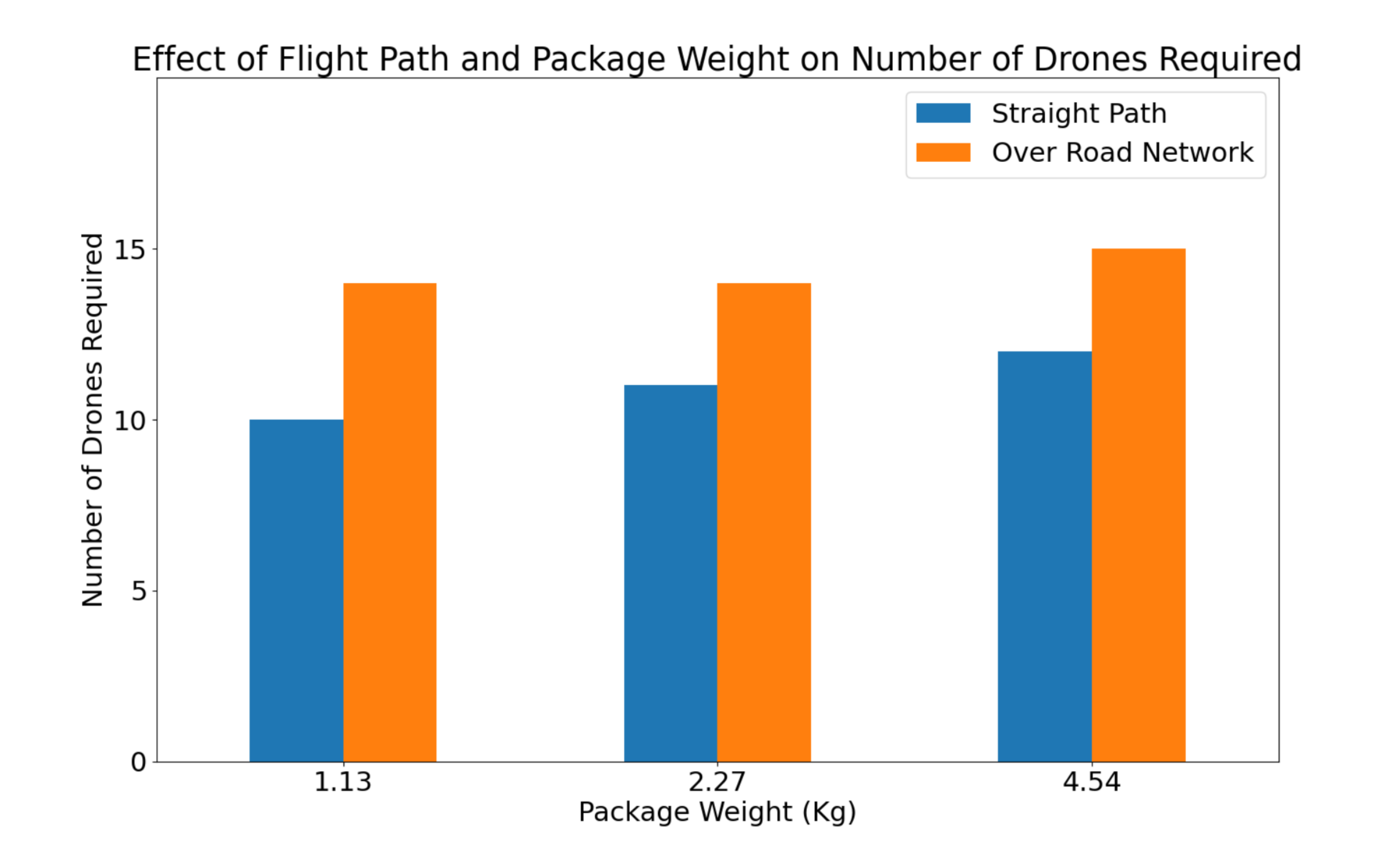}}\caption{\label{fig:flight_path_number_of_drones}Effect of flight path and
package weight on the required number of drones.}
\par\end{centering}
\end{figure}

Figures \ref{fig:flight_path_speed_30_number_of_drones} and \ref{fig:flight_path_speed_15mph_number_of_drones}
demonstrate the effect of flight path and package weight on the required
fleet size (i.e., required number of drones) to deliver packages to
all customer locations with drone speeds of 13.41 m/s and 6.71 m/s,
respectively. We see that a larger number of drones are needed as
the drones fly over the road networks as compared to the straight
path at both speed levels. This is because, as compared to the straight
path, while flying over the road networks, drones travel much longer
distances to deliver packages to customer locations, increasing the
travel time. With this increased travel time over the road networks
for each delivery, it is not possible to meet the pickup time window
with the same number of drones as in flying in the straight path with
shorter distance and travel time. Therefore, a larger number of drones
are needed to deliver the packages meeting the time window requirement
while flying over the road networks. When the DJI Matrice 600 Pro
drones are flying over the road networks at 13.41 m/s speed, the percentage
increase in the required number of drones across the three different
package weights are 22.2\%, 22.2\%, and 33.3\% as compared to flying
in the straight path. For the 6.71 m/s drone speed, the percentage
increase in the required number of drones are 40.0\%, 27.3\%, and
25.0\% as compared to flying in the straight path corresponding to
the three different package weights. However, it is evident from Figures
\ref{fig:flight_path_speed_30_number_of_drones} and \ref{fig:flight_path_speed_15mph_number_of_drones}
that the package weight does not have much noticeable impact on the
number of drones as the package weight does not affect drone speed
and travel time.

\begin{figure}[h]
\begin{centering}
\subfloat[\label{fig:flight_path_energy-consumption_30mph}Total energy consumption
for drone speed of 13.41 m/s]{\begin{centering}
\includegraphics[width=8cm,height=8cm,keepaspectratio]{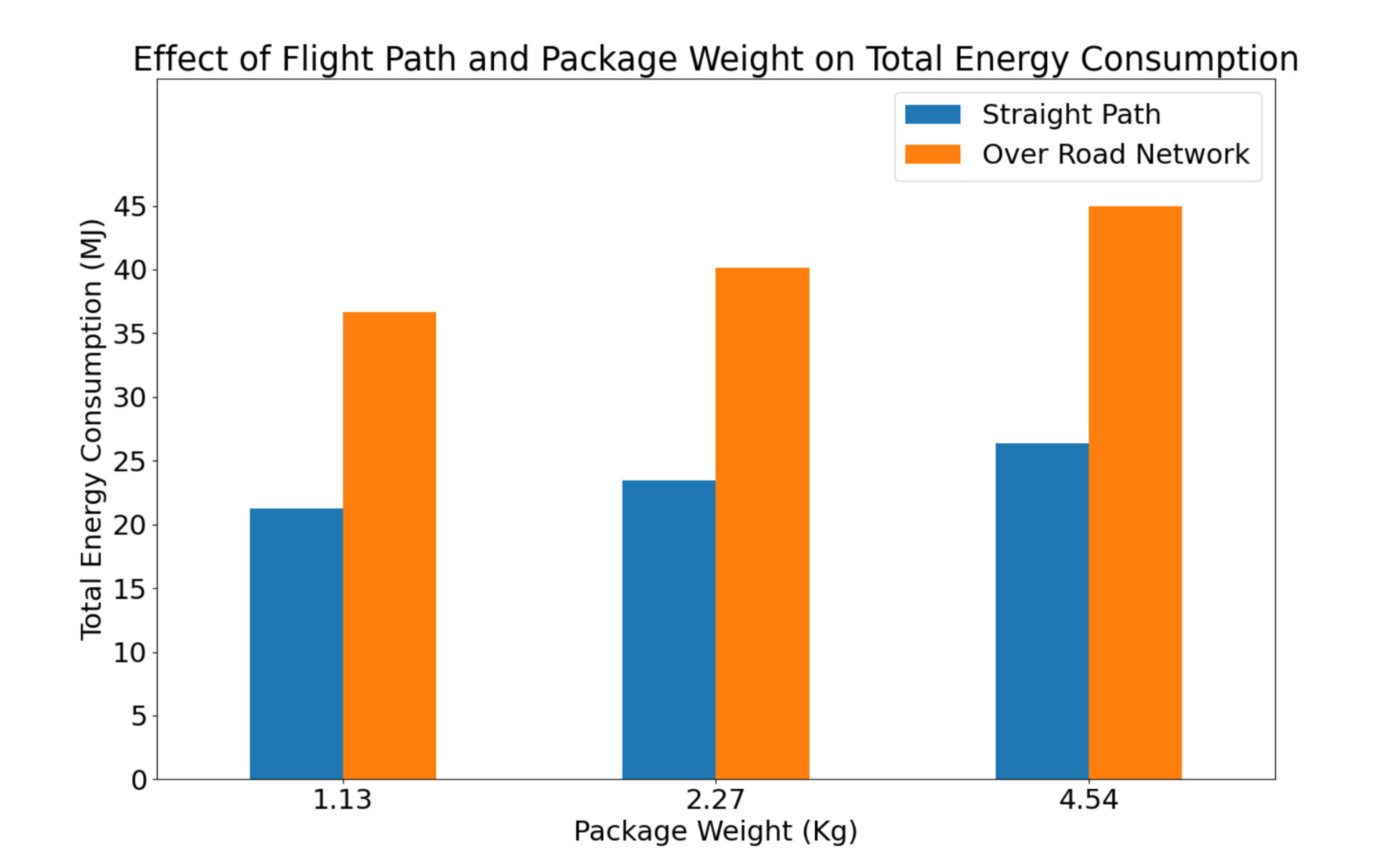}
\par\end{centering}

}~~\subfloat[\label{fig:flight_path_energy-consumption_15mph}Total energy consumption
for drone speed of 6.71 m/s ]{\begin{centering}
\includegraphics[width=8cm,height=8cm,keepaspectratio]{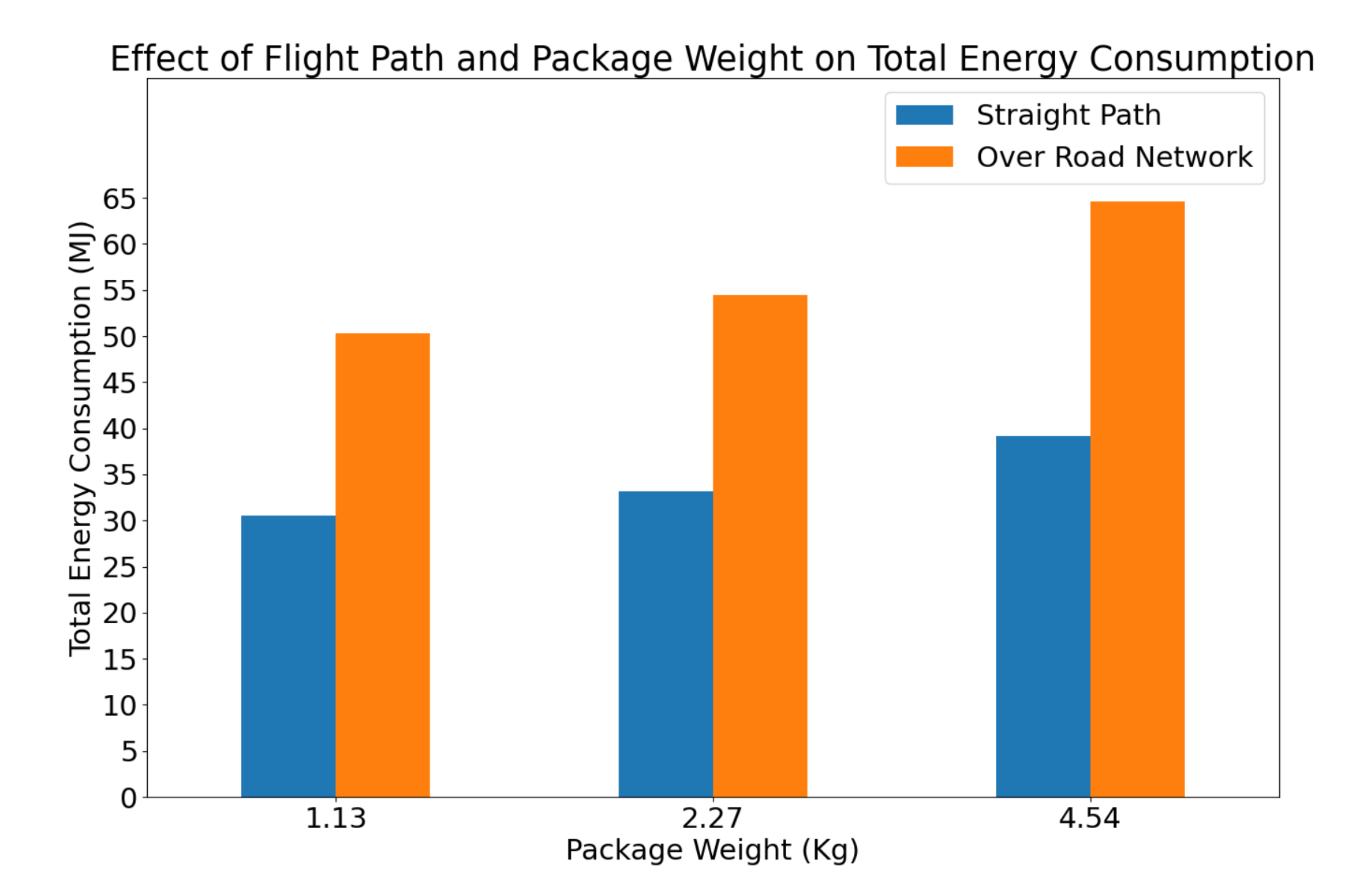}
\par\end{centering}

}\caption{\label{fig:Flight_path_energy_consumption}Effect of flight path and
package weight on total energy consumption.}
\par\end{centering}
\end{figure}

Figures \ref{fig:flight_path_energy-consumption_30mph} and \ref{fig:flight_path_energy-consumption_15mph}
demonstrate the effect of flight path and package weight on the total
energy consumption in delivering all customer orders at the drone
speeds of 13.41 m/s and 6.71 m/s, respectively. We see that the total
energy consumption increases substantially as the drones fly over
the road networks as compared to flying in a straight path at both
speed levels across all three package weights. This is because flying
over the road networks requires the drones to travel a longer distance
and stay in the air for a longer amount of time, resulting in an increased
energy consumption. We have observed that the total energy consumption
increases by 72.22\%, 71.24\%, and 70.54\% corresponding to the package
weights of 1.13 Kg, 2.27 Kg, and 4.54 Kg, respectively, as the drones
fly over the road networks than in the straight paths at 13.41 m/s
speed. The percentage increase in the total energy consumption due
to flying over the road networks at 6.71 m/s speed are 64.85\%, 64.14\%,
and 64.87\% across three different package weights. The average energy
consumption per delivery also demonstrate the same trend and magnitude
(see Figure \ref{fig:Flight_path_avg_energy_consumption} in Appendix
\ref{sec:Effect_flight_path_average_energy_consumption}). It is evident
from Figures \ref{fig:flight_path_energy-consumption_30mph} and \ref{fig:flight_path_energy-consumption_15mph}
that the energy consumption increases as the drones carry larger package
weight at both speed levels, which is very intuitive.

\begin{figure}[h]
\begin{centering}
\subfloat[\label{fig:flight_path_battery_replacements_30mph}Required number
of battery replacements for drone speed of 13.41 m/s.]{\begin{centering}
\includegraphics[width=8cm,height=8cm,keepaspectratio]{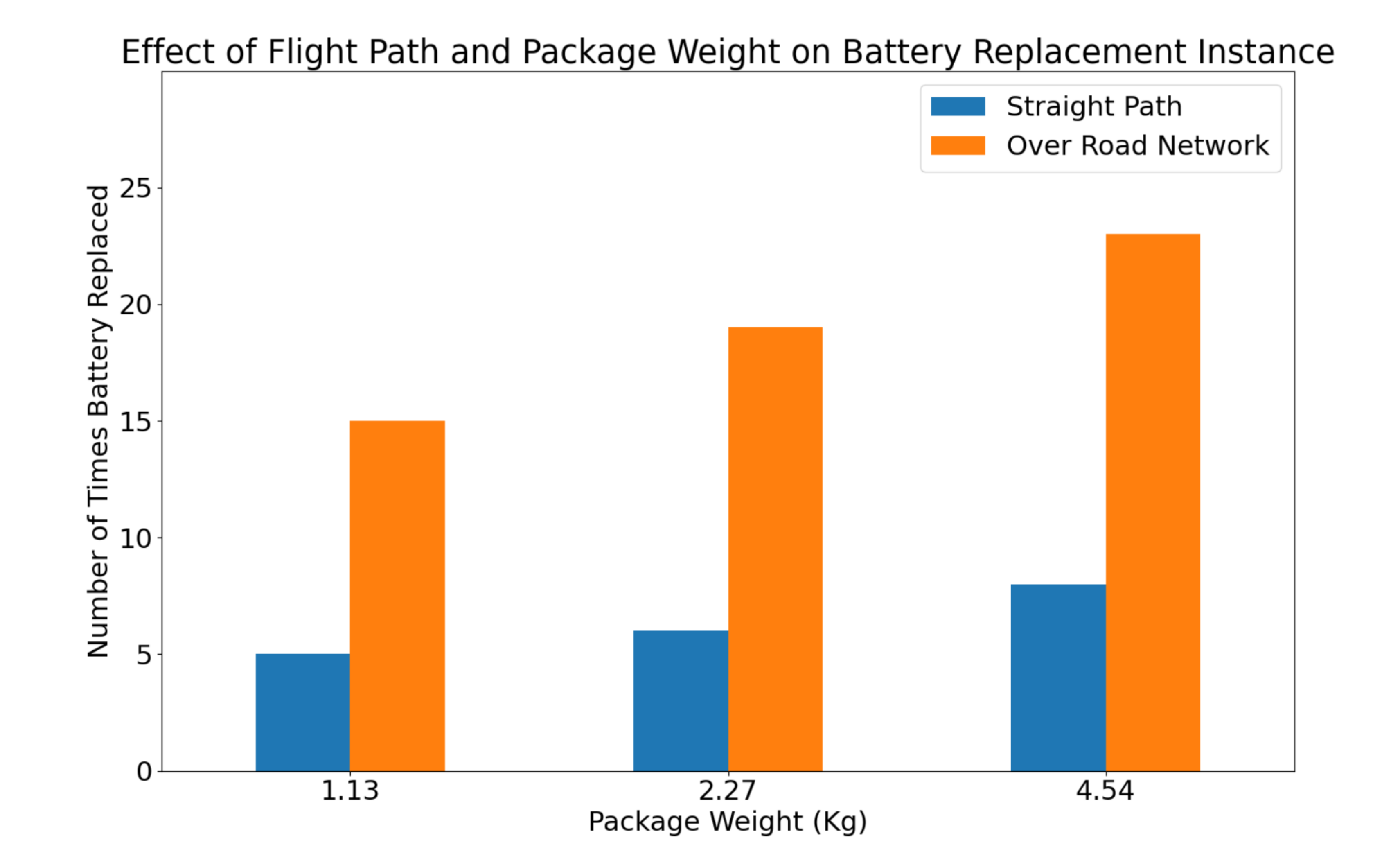}
\par\end{centering}

}~~\subfloat[\label{fig:flight_path_battery_replacements_15mph}Required number
of battery replacements for drone speed of 6.71 m/s.]{\begin{centering}
\includegraphics[width=8cm,height=8cm,keepaspectratio]{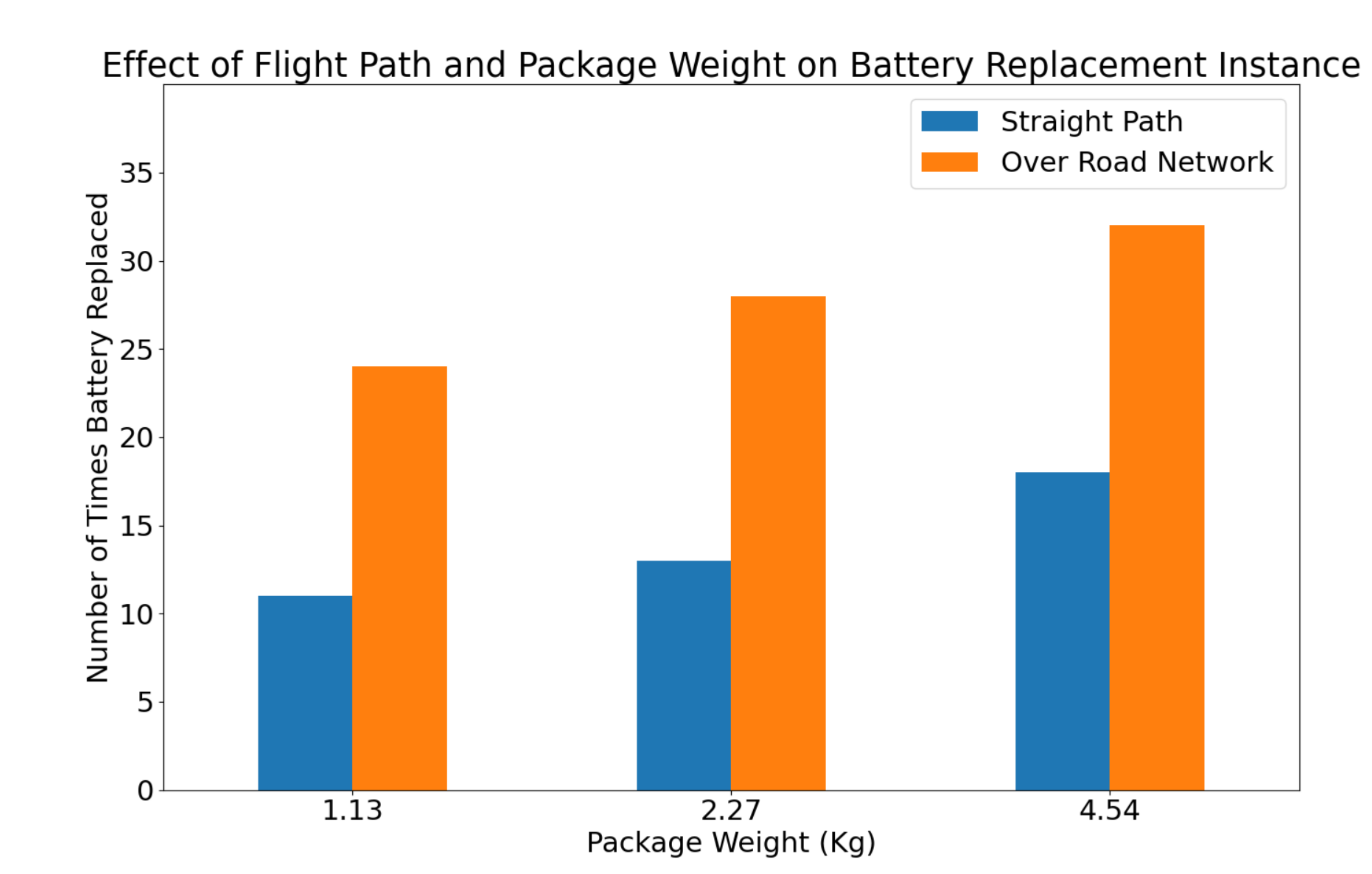}
\par\end{centering}
}\caption{Effect of flight path and package weight on the required number of
battery replacements.}
\par\end{centering}
\end{figure}

Figures \ref{fig:flight_path_battery_replacements_30mph} and \ref{fig:flight_path_battery_replacements_15mph}
demonstrate that the required number of battery replacements increases
substantially as the drones fly over the road networks as compared
to a straight path for both speed levels. For instance, the required
number of battery replacements increases by 200\% as the drones fly
over the road networks compared to flying in a straight path at 13.41
m/s speed and all customer orders weigh 1.13 Kg. This is due to the
increased energy consumption while flying over the road networks as
explained earlier. As more energy is consumed in delivering each customer
order, the minimum required battery energy is reached more frequently,
resulting in a larger number of battery replacements.

\subsection{Effect of Mixed Drone Fleet on Energy Consumption and Required Number
of Drones}

We observed that due to its small size and limited battery capacity,
the delivery range and package weight carrying capacity of the Tarot
650 drone (rotary quadcopter) is much smaller compared to the DJI
Matrice 600 Pro (rotary hexacopter). But our drone flight test data
shows that the Tarot 650 drone is more energy-efficient than the DJI
Matrice 600 Pro drone while carrying the same package weight for the
same distance. Therefore, our assumption is that business owners can
substantially reduce energy consumption by using a mixed fleet of
drones consisting of these both types of drones such that the Tarot
650 drones are used to deliver smaller packages to shorter distances
whereas larger packages are delivered longer distances using the DJI
Matrice 600 Pro drones.

Using the modeling approach presented in Section \ref{subsec:Modeling_approach_mixed_fleet},
we evaluated the effect of using both Tarot and DJI drones in a mixed
fleet on the number of locations served, energy consumption, and the
required number of drones of each type. As the maximum weight carrying
capacity of the Tarot 650 drone is 1.13 Kg, we used the package weight
of 1.13 Kg for both types of drones in the analysis. Table \ref{tab:Mixed_fleet_delivery_locations_served}
demonstrates that 113 of the 126 locations are served using the Tarot
drones, whereas only 13 locations are served using DJI drones when
all drones are operated at 13.41 m/s speed. We see that the Tarot
drones contribute 61.95\% to the total energy consumption in delivering
packages to 113 locations, whereas the DJI drones contribute 38.05\%
despite delivering packages to only 13 locations. DJI drones are only
used to deliver packages to locations outside the delivery range of
Tarot. We see a similar pattern in the number of locations served
and total energy consumption by the mixed fleet of drones when the
drones are operated at 6.71 m/s speed, as shown in Table \ref{tab:Mixed_fleet_delivery_locations_served},
which also demonstrates the optimal number of drones required for
each type of drone in the mixed fleet.

\begin{table}[h]
\caption{\label{tab:Mixed_fleet_delivery_locations_served}Number of locations
served, percentage of total energy consumption, and the required number
of drones for different drone types in a mixed fleet.}

\begin{centering}
{\small{}}%
\begin{tabular}[t]{>{\centering}m{1.5cm}>{\centering}p{2cm}>{\centering}p{2cm}>{\centering}p{2cm}|>{\centering}p{2cm}>{\centering}p{2cm}>{\centering}p{2cm}}
\hline 
\multirow{2}{1.5cm}{\textbf{Drone Speed (m/s)}} & \multicolumn{3}{c|}{\textbf{DJI Matrice 600 Pro}} & \multicolumn{3}{c}{\textbf{Tarot 650}}\tabularnewline
\cline{2-7} \cline{3-7} \cline{4-7} \cline{5-7} \cline{6-7} \cline{7-7} 
 & \multirow{1}{2cm}{\centering{}\textbf{Number of Locations Served}} & \textbf{\% of Total Energy Consumption} & \textbf{Required Numbers of Drones} & \centering{}\textbf{Number of Locations Served} & \textbf{\% of Total Energy Consumption} & \textbf{Required Numbers of Drones}\tabularnewline
\hline 
\hline 
\centering{}\textbf{13.41} & \centering{}13 & 38.05  & 5 & \centering{}113 & \centering{}61.95  & 30\tabularnewline
\centering{}\textbf{6.71} & \centering{}21 & \centering{}53.74  & 9 & 77 & 46.26  & 21\tabularnewline
\hline 
\end{tabular}{\small\par}
\par\end{centering}
\end{table}

\begin{figure}
\begin{centering}
\includegraphics[width=12cm,height=12cm,keepaspectratio]{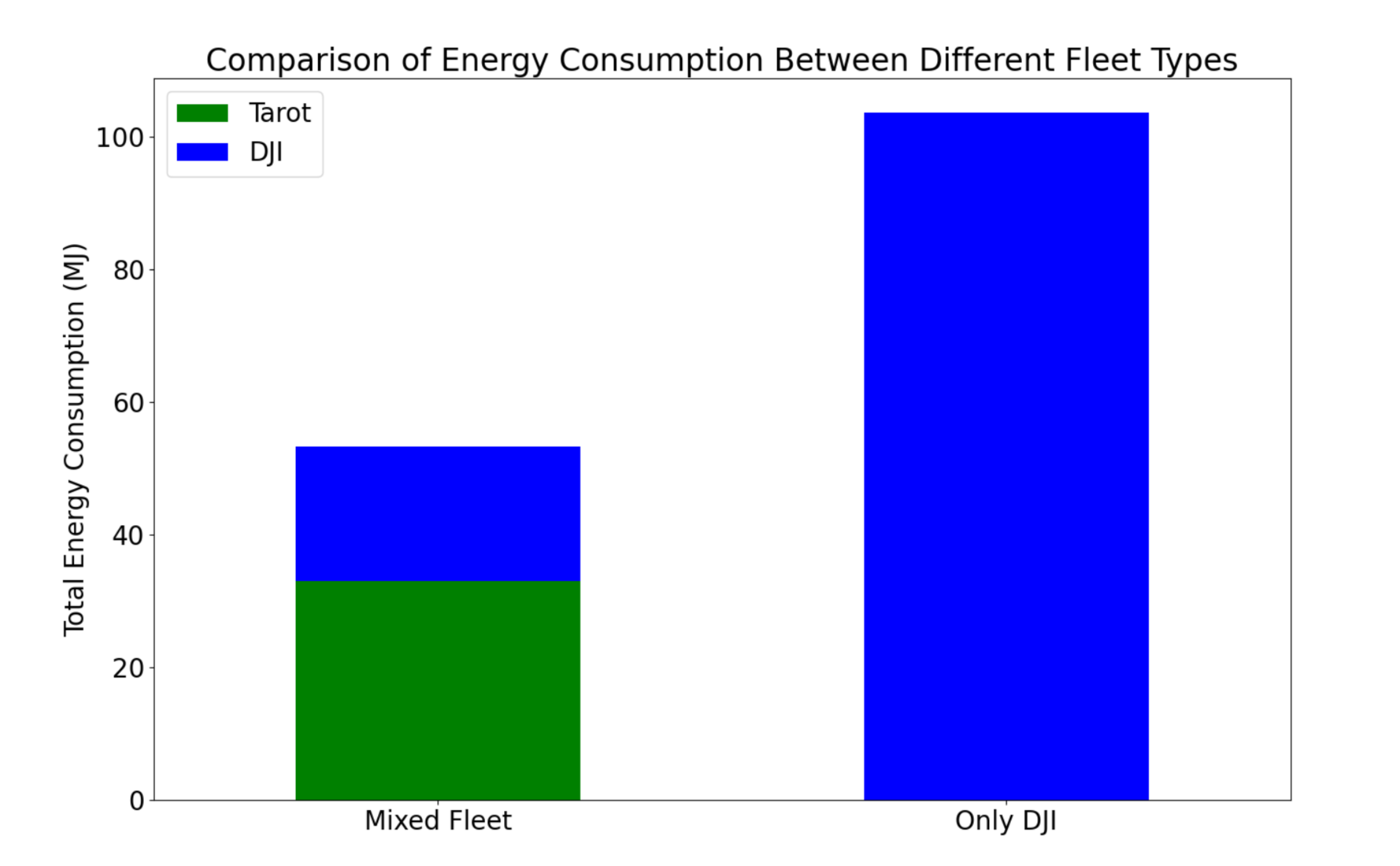}\caption{\label{fig:energy_consumption_comparison_mixed_fleet}Comparison of
total energy consumption between different fleet types.}
\par\end{centering}
\end{figure}

We analyzed how much does total energy consumption decreased due to
using a mixed fleet of drones compared to using a homogeneous fleet
with only DJI drones. Figure \ref{fig:energy_consumption_comparison_mixed_fleet}
demonstrates the total energy consumption by the mixed fleet of drones
and the homogeneous fleet with only DJI drones when all drones are
operated at 13.41 m/s speed. We see that using a mixed fleet with
both Tarot and DJI drones substantially reduces the energy consumption
in delivering customer orders compared to the homogeneous fleet with
only DJI drones. For instance, when the drones are operated at 13.41
m/s speed, using a mixed fleet of drones reduces the total energy
consumption by 48.52\% compared to the homogeneous fleet with only
DJI drones.

\section{Conclusion}

We studied a drone deployment optimization problem for direct delivery
of goods from a retail business location that acts as the central
depot to customer locations, where each drone visits a single customer
location on each trip and returns to the depot before visiting the
next customer. In this problem, the goal of the business owner is
to optimally route the drones that minimizes the total energy consumption,
the required fleet size (i.e., required number of drones), and the
required number of battery replacements in delivering customer orders
maintaining a specified pickup time window. We proposed a new mathematical
optimization approach to solve the problem for both homogeneous and
mixed fleets of drones. The mathematical model incorporates the specified
time window, factors affecting the drone energy consumption, and the
drone operating parameters (e.g., drone speed, package weight, flight
path, the minimum required battery energy, and fleet type) and optimally
decides the routing of drones so that the minimum number of drones
are used, the number of times drone battery is replaced is minimum
and the total energy consumption is minimum in delivering customer
orders. We also proposed two valid inequalities to improve the computational
efficiency of the mathematical optimization models.

We conducted numerical experiments using real drone flight test data
to provide new insights into the effect of different drone operating
conditions on the performance of drone-based delivery. Results demonstrate
that the total and average energy consumption increases as the drones
fly at a slower speed due to staying in the air for a longer time.
For instance, the total energy consumption for DJI drones (i.e., rotary
hexacopter) increases by 47.96\% as the drone speed reduces from 13.41
m/s to 6.71 m/s when all customer orders weigh 1.13 Kg. Energy consumption
increases substantially as the drones fly over the road network as
compared to a straight path; the total energy consumption for DJI
drones increases by 72.22\% as the drones fly over the road networks
compared to the straight path with all customers weighing 1.13 Kg.
Drone energy consumption also increases with package weight and package
dropping delivery method that requires a longer hover duration. Results
demonstrate that using a mixed fleet of drones (i.e., a fleet containing
both rotary hexacopters and quadcopters) reduces the total energy
consumption by 48.52\% compared to the homogeneous fleet with only
hexacopter drones when all drones are operated at 13.41 m/s speed.
The delivery range of drones decreases as the package weight increases,
drone speed decreases, and the minimum required battery energy increases.

The required fleet size increases as the pickup time window of the
customer orders decreases, as a higher number of drones are needed
to pick up the products within a shorter time window. The number of
drones required increases by 48\% as the time window reduces from
35 minutes to 1 minute. The required number of drones increases as
the drone speed decreases, as a higher number of drones are required
to maintain the same pickup time window when the drones are flying
at a slower speed. A higher number of drones are needed when the drones
fly over the road networks compared to flying in a straight path;
the required number of drones increases by 22.2\% as the drones fly
over the road networks as compared to flying in a straight path.

The required number of battery replacements increases as the package
weight and minimum required battery energy increases and speed decreases.
As the drones fly over the road networks, the drone battery needs
to be replaced a substantially larger number of times as compared
to flying in a straight path. This is mainly because of the higher
energy consumption due to traveling long distances over the road network.
The required number of battery replacements increases by 200\% as
the drones fly over the road network compared to flying in a straight
path when all customer orders weigh 1.13 Kg and drones fly at 13.41
m/s speed. Therefore, business owners need to have a substantially
larger number of additional drone batteries when their drones need
to be flown over the road networks as compared to flying in a straight
path.

The managerial insights drawn from the numerical results based on
real drone flight tests and delivery data will help retail business
owners better understand the potential of drones under different conditions.
Another practical implication of our research is that retail business
owners can use our proposed optimization method to make investment
and operating decisions in using their drones to deliver customer
orders. Specifically, a retail business owner can use our optimization
method with their drone energy profile and delivery data to decide
on the minimum required fleet size (i.e., the required number of drones)
and the number of drones in each type for a mixed fleet under different
conditions (e.g., time window requirement, flight path, drone speed).
The business owners can compute the number of additional batteries
they need to have to operate the drone fleet depending on their drone
operating conditions. Additionally, business owners can evaluate how
different drone operating parameters specific to their business and
products affect the total energy consumption.

\section*{Acknowledgments}

We thank the U.S. Department of Energy Vehicle Technology Office (VTO)
for funding and supporting this work. This work is supported by the
U.S. Department of Energy under Department of Energy Idaho Operations
Office Contract No. DE-AC07-05ID14517.

\bibliographystyle{unsrt}
\bibliography{Manuacript_scenario1}

\newpage

\section*{Appendices}

\renewcommand\thefigure{\thesection.\arabic{figure}}
\setcounter{figure}{0}
\renewcommand\thetable{\thesection.\arabic{table}}
\setcounter{table}{0}

\appendix

\section{\label{sec:Flowchart_mixed_fleet}Flowchart of the Mixed Fleet of
Drones Modeling}

\begin{figure}[h]

\begin{centering}
\includegraphics[width=18cm,height=18cm,keepaspectratio]{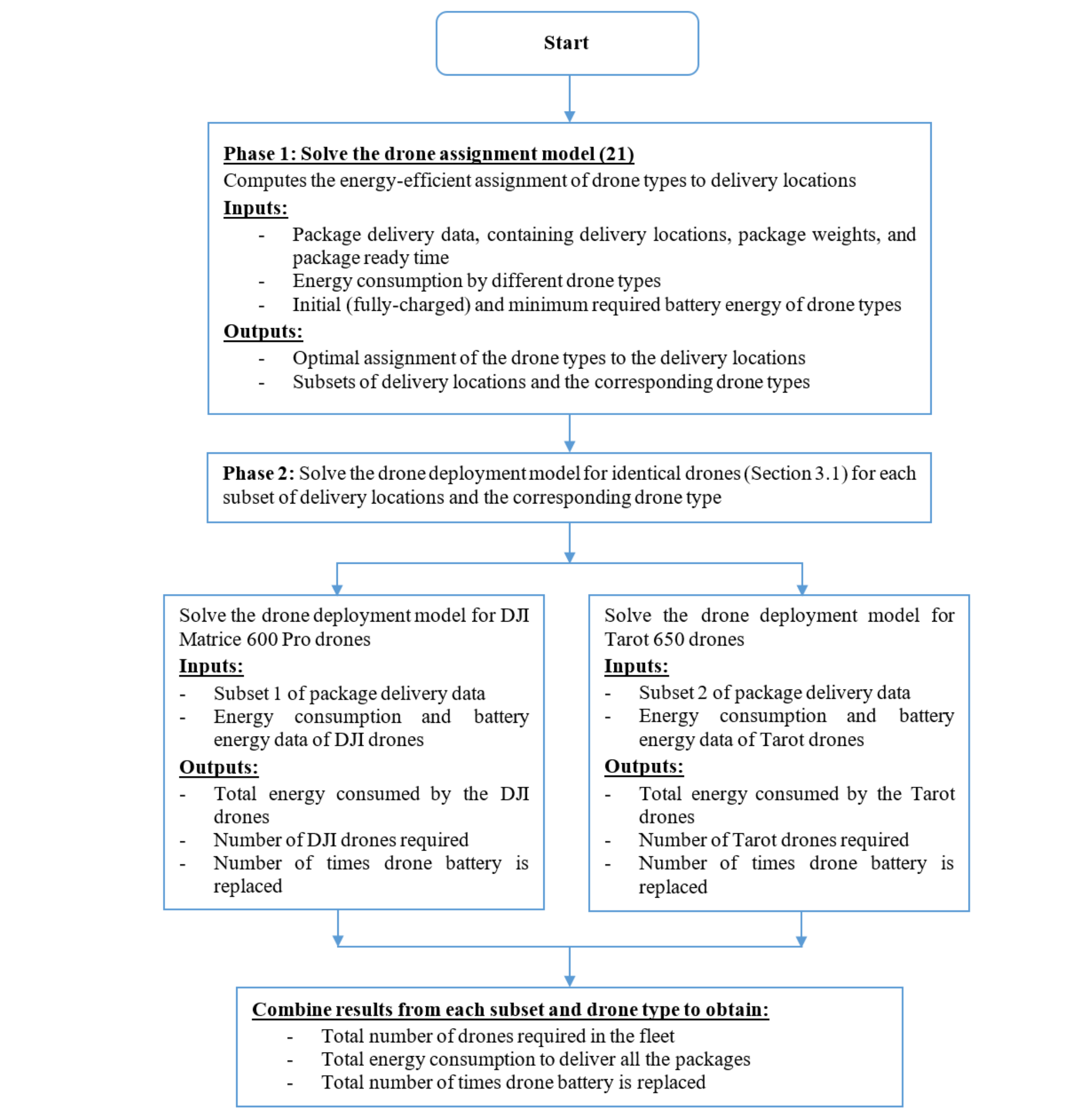}
\par\end{centering}
\caption{\label{fig:Two-phase-modeling-approach}Two-phase modeling approach
for mixed fleet of drones.}

\end{figure}

Figure \ref{fig:Two-phase-modeling-approach} demonstrates the two-phase
modeling approach proposed to solve the drone deployment optimization
problem for a mixed fleet of drones. 

\section{Energy Consumption and Flight Time Data\label{sec:Energy-Consumption_time_segments}}

Table \ref{tab:Energy-consumption_segments-DJI} presents the average
power consumption data for different flight segments of our DJI Matrice
600 Pro drone at two different speed levels and four different package
weights. 

\begin{table}[h]
\caption{\label{tab:Energy-consumption_segments-DJI}Power consumption of the
DJI Matrice 600 Pro drone for the flight segments.}

\begin{centering}
{\small{}}%
\begin{tabular}{>{\raggedright}m{3cm}>{\raggedright}m{3cm}>{\centering}p{2cm}>{\centering}p{2cm}>{\centering}p{2cm}>{\centering}p{2cm}}
\hline 
\multirow{2}{3cm}{\centering{}\textbf{Drone Speed (m/s)}} & \multirow{2}{3cm}{\centering{}\textbf{Package Weight (Kg)}} & \multicolumn{4}{c}{\textbf{Average Power (Joule per Second) in Flight Segments}}\tabularnewline
\cline{3-6} \cline{4-6} \cline{5-6} \cline{6-6} 
 &  & \centering{}\textbf{Ascend} & \centering{}\textbf{Descend} & \centering{}\textbf{Forward Flight} & \textbf{Hover}\tabularnewline
\hline 
\hline 
\multirow{4}{3cm}{\centering{}13.41} & \centering{}0 & \centering{}1351.4456 & \centering{}1023.8680 & \centering{}1186.5048 & 1039.2542\tabularnewline
 & \centering{}1.13 & \centering{}1487.3006 & \centering{}1104.4719 & \centering{}1479.2276 & 1211.6308\tabularnewline
 & \centering{}2.27 & 1746.2067 & 1422.7263 & 1718.7145 & 1406.6182\tabularnewline
 & \centering{}4.54 & 2233.5910 & 1738.2538 & 2044.0478 & 1759.3900\tabularnewline
\hline 
\multirow{4}{3cm}{\centering{}6.71} & \centering{}0 & \centering{}1351.4456 & \centering{}1023.8680 & \centering{}1052.5789 & 1039.2542\tabularnewline
 & \centering{}1.13 & \centering{}1487.3006 & \centering{}1104.4719 & 1232.8998 & 1211.6308\tabularnewline
 & \centering{}2.27 & 1746.2067 & 1422.7263 & 1401.8190 & 1406.6182\tabularnewline
 & \centering{}4.54 & 2233.5910 & 1738.2538 & \centering{}1878.1490 & 1759.3900\tabularnewline
\hline 
\end{tabular}{\small\par}
\par\end{centering}
\end{table}

Table \ref{tab:Time-duration_segments-DJI} presents the flight time
of different flight segments of our DJI Matrice 600 Pro drone test
at two different speed levels. Flight durations are the same for different
package weights at the same speed level. 

\begin{table}
\caption{\label{tab:Time-duration_segments-DJI}Flight time of the DJI Matrice
600 Pro drone for the flight segments.}

\begin{centering}
{\small{}}%
\begin{tabular}{>{\raggedright}m{3cm}>{\centering}p{2cm}>{\centering}p{2cm}>{\centering}p{2cm}>{\centering}p{2cm}>{\centering}p{2cm}}
\hline 
\multirow{3}{3cm}{\centering{}\textbf{Drone Speed (m/s)}} & \multicolumn{5}{c}{\textbf{Duration (Seconds) of Flight Segments}}\tabularnewline
\cline{2-6} \cline{3-6} \cline{4-6} \cline{5-6} \cline{6-6} 
 & \multirow{2}{2cm}{\centering{}\textbf{Ascend}} & \multirow{2}{2cm}{\centering{}\textbf{Descend}} & \multicolumn{2}{c}{\textbf{Forward Flight }} & \multirow{2}{2cm}{\centering{}\textbf{Hover}}\tabularnewline
\cline{4-5} \cline{5-5} 
 &  &  & \centering{}\textbf{Time/mile} & \centering{}\textbf{Time/kilometer} & \tabularnewline
\hline 
\hline 
\multirow{1}{3cm}{\centering{}13.41} & \centering{}24.6 & \centering{}41.80 & \centering{}125.0 & 78.125 & 5.0\tabularnewline
\hline 
\multirow{1}{3cm}{\centering{}6.71} & \centering{}24.6 & \centering{}41.80 & \centering{}245.0 & 153.125 & 5.0\tabularnewline
\hline 
\end{tabular}{\small\par}
\par\end{centering}
\end{table}

Table \ref{tab:Energy-consumption_segments-TAROT} presents the average
power consumption data for different flight segments of our Tarot
650 drone test at two different speed levels and two different package
weights. 

\begin{table}[h]
\caption{\label{tab:Energy-consumption_segments-TAROT}Power consumption of
the Tarot 650 drone for the flight segments.}

\centering{}{\small{}}%
\begin{tabular}{>{\raggedright}m{3cm}>{\raggedright}m{3cm}>{\centering}p{2cm}>{\centering}p{2cm}>{\centering}p{2cm}>{\centering}p{2cm}}
\hline 
\multirow{2}{3cm}{\centering{}\textbf{Drone Speed (m/s)}} & \multirow{2}{3cm}{\centering{}\textbf{Package Weight (Kg)}} & \multicolumn{4}{c}{\textbf{Average Power (Joule per Second) in Flight Segments}}\tabularnewline
\cline{3-6} \cline{4-6} \cline{5-6} \cline{6-6} 
 &  & \centering{}\textbf{Ascend} & \centering{}\textbf{Descend} & \centering{}\textbf{Forward Flight} & \textbf{Hover}\tabularnewline
\hline 
\hline 
\multirow{2}{3cm}{\centering{}13.41} & \centering{}0 & \centering{}419.9395 & \centering{}370.3456 & \centering{}389.3249 & \centering{}369.2475\tabularnewline
 & \centering{}1.13 & \centering{}638.3887 & \centering{}589.3817 & \centering{}596.6649 & 564.6385\tabularnewline
\hline 
\multirow{2}{3cm}{\centering{}6.71} & \centering{}0 & \centering{}419.9395 & \centering{}370.3456 & \centering{}381.6552 & \centering{}369.2475\tabularnewline
 & \centering{}1.13 & \centering{}638.3887 & \centering{}589.3817 & 580.4716  & 564.6385\tabularnewline
\hline 
\end{tabular}{\small\par}
\end{table}

Table \ref{tab:Time-duration_segments-TAROT} presents the flight
time of different flight segments of our Tarot 650 drone test at two
different speed levels. Flight durations are the same for different
package weights at the same speed level. 

\begin{table}
\caption{\label{tab:Time-duration_segments-TAROT}Flight time of the Tarot
650 drone for the flight segments.}

\centering{}{\small{}}%
\begin{tabular}{>{\raggedright}m{3cm}>{\centering}p{2cm}>{\centering}p{2cm}>{\centering}p{2cm}>{\centering}p{2cm}>{\centering}p{2cm}}
\hline 
\multirow{3}{3cm}{\centering{}\textbf{Drone Speed (m/s)}} & \multicolumn{5}{c}{\textbf{Duration (Seconds) of Flight Segments}}\tabularnewline
\cline{2-6} \cline{3-6} \cline{4-6} \cline{5-6} \cline{6-6} 
 & \multirow{2}{2cm}{\centering{}\textbf{Ascend}} & \multirow{2}{2cm}{\centering{}\textbf{Descend}} & \multicolumn{2}{c}{\textbf{Forward Flight}} & \multirow{2}{2cm}{\centering{}\textbf{Hover}}\tabularnewline
\cline{4-5} \cline{5-5} 
 &  &  & \centering{}\textbf{Time/mile} & \centering{}\textbf{Time/kilometer} & \tabularnewline
\hline 
\hline 
\multirow{1}{3cm}{\centering{}13.41} & \centering{}32.56 & \centering{}38.31 & \centering{}130.0 & 80.78 & 5.0\tabularnewline
\hline 
\multirow{1}{3cm}{\centering{}6.71} & \centering{}32.56 & \centering{}38.31 & \centering{}250.0 & 155.34 & 5.0\tabularnewline
\hline 
\end{tabular}{\small\par}
\end{table}

\section{Effect of Drone Flight Path on Average Energy Consumption\label{sec:Effect_flight_path_average_energy_consumption} }

\begin{figure}[h]
\centering{}\subfloat[\label{fig:flight_path_average_energy-consumption_30mph}Average energy
consumption for drone speed of 13.41 m/s.]{\begin{centering}
\includegraphics[width=8cm,height=8cm,keepaspectratio]{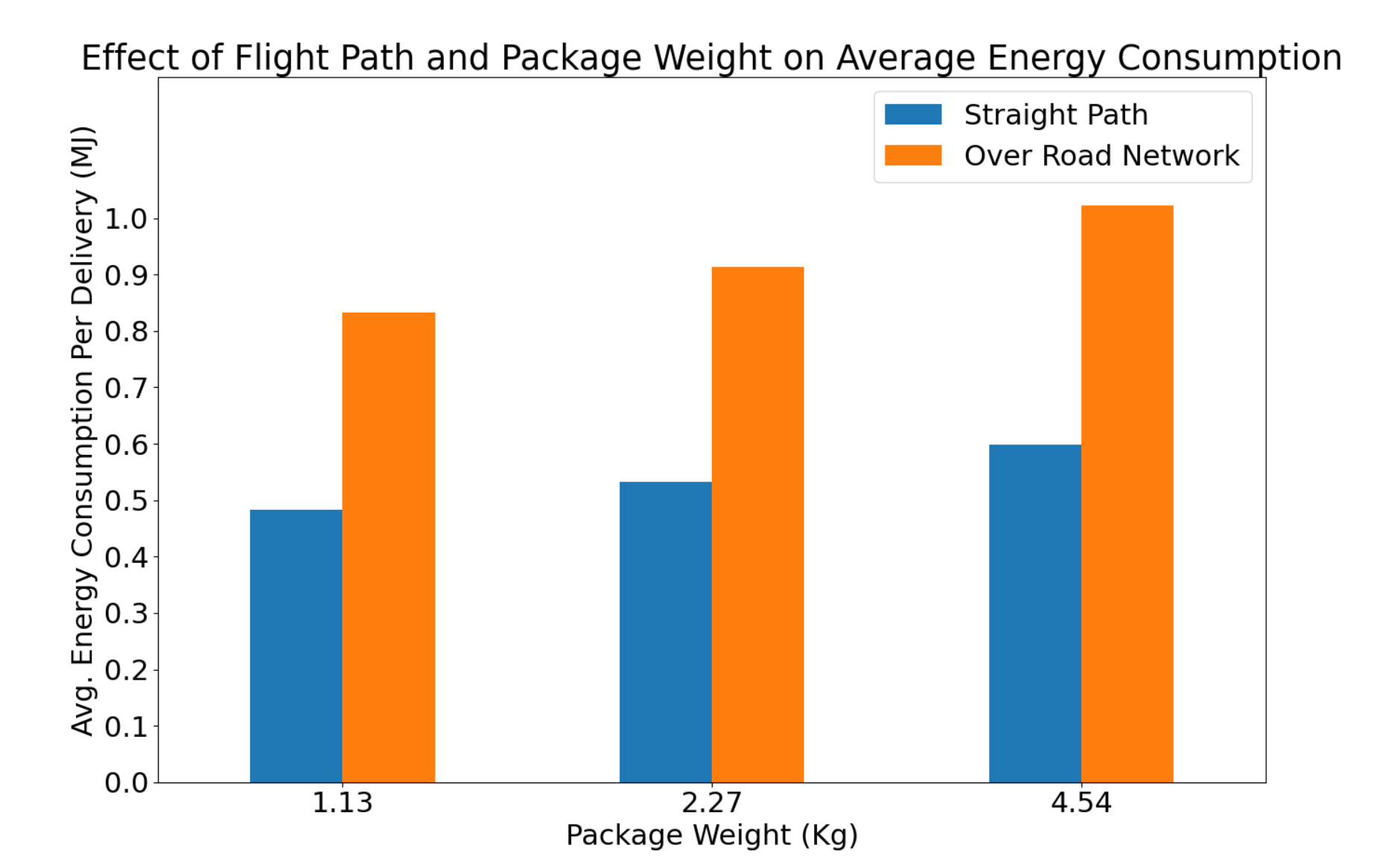}
\par\end{centering}
}~~\subfloat[\label{fig:flight_path_average_energy-consumption_15mph}Average energy
consumption for drone speed of 6.71 m/s. ]{\begin{centering}
\includegraphics[width=8cm,height=8cm,keepaspectratio]{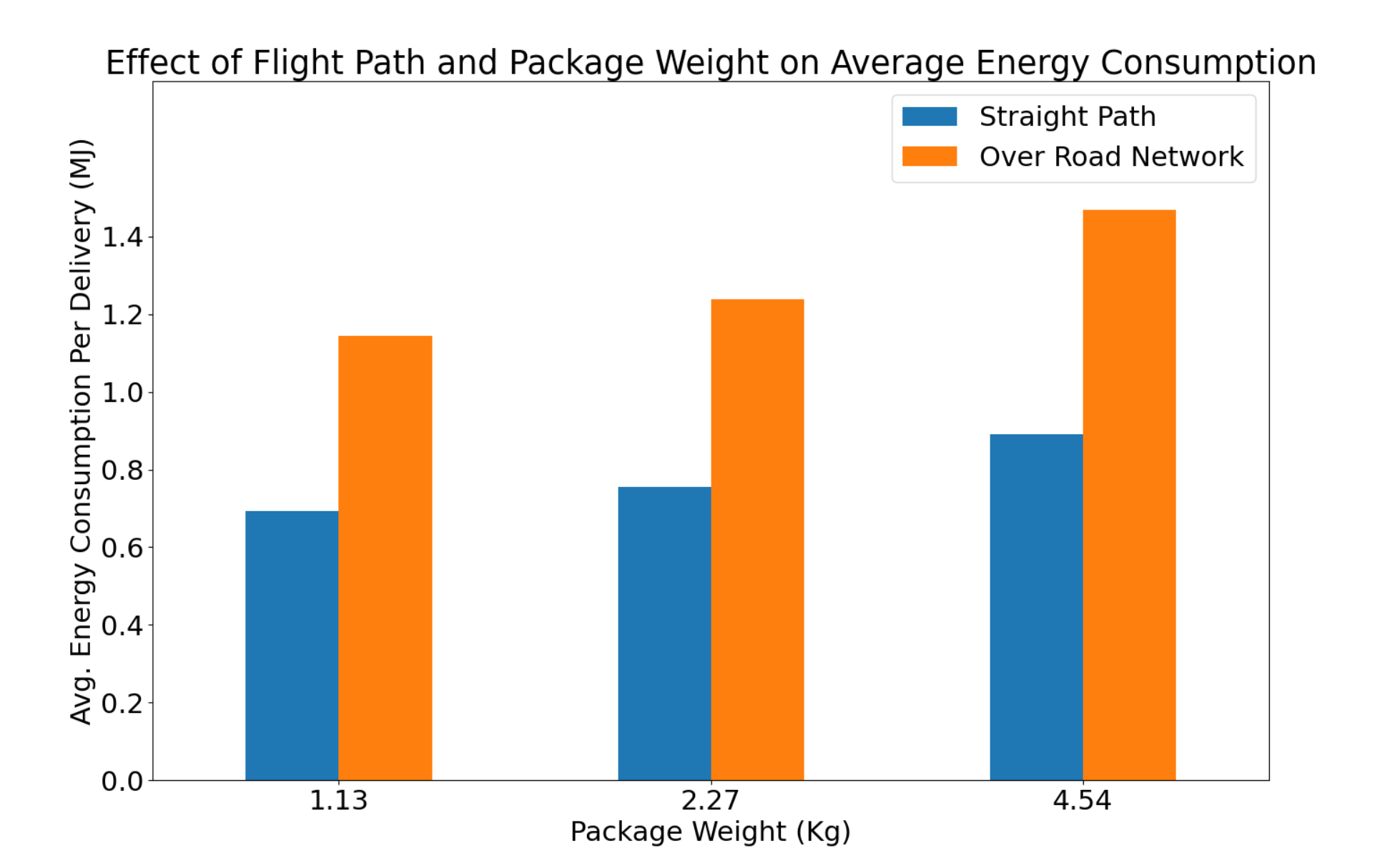}
\par\end{centering}
}\caption{\label{fig:Flight_path_avg_energy_consumption}Effect of flight path
and package weight on the average energy consumption per delivery.}
\end{figure}

Figures \ref{fig:flight_path_average_energy-consumption_30mph} and
\ref{fig:flight_path_average_energy-consumption_15mph} demonstrate
the effect of flight path and package weight on the average energy
consumption per delivery of the DJI Matrice 600 Pro drone at drone
speeds of 13.41 m/s and 6.71 m/s, respectively. The average energy
consumption increases by 72.22\%, 71.25\%, and 70.54\% corresponding
to the package weights of 1.13 Kg, 2.27 Kg, and 4.54 Kg, respectively,
as the drones fly over the road networks than the straight paths at
13.41 m/s speed. The percentage increase in the average energy consumption
due to flying over the road networks at 6.71 m/s speed are 64.48\%,
64.13\%, and 64.87\% across three different package weights. 

\end{document}